\documentclass[hyperref]{ctexart}
\usepackage{geometry} 
\usepackage{helvet}
\usepackage{amsmath, amsfonts, amssymb}
\usepackage[english]{babel}

\usepackage{url}      
\usepackage{bm}      
\usepackage{multirow}
\usepackage{booktabs}
\usepackage{cite}
\usepackage{lipsum}
\usepackage{graphicx}
\usepackage{graphics}
\usepackage{subfigure}
\usepackage{epstopdf}
\usepackage{float}
\usepackage{indentfirst}
\usepackage{algorithm}
\usepackage{algorithmicx}  
\usepackage{algpseudocode}

\usepackage{fancyhdr} 
\usepackage{hyperref} 
\hypersetup{colorlinks, bookmarks, unicode} 

\numberwithin{equation}{section}

\newtheorem{Theorem}{Theorem}[section]
\newtheorem{Lemma}{Lemma}[section]

\newtheorem{Property}{Property}[section]
\newtheorem{Monotonicity Criterion}{Monotonicity Criterion}

\newtheorem{Remark}{Remark}[section]

\newtheorem{Definition}{Definition}[section]
\newtheorem{Proposition}{Proposition}[section]
\newtheorem{Proof}{Proof}[section]

\title{\textbf{A well-balanced scheme for Euler equations with singular sources}}
\author{\sffamily Changsheng Yu, \sffamily Tiegang Liu, \sffamily Chengliang Feng}

\begin{document}
	\maketitle
	
\begin{abstract}
	\begin{abstract}
		Numerical methods for the Euler equations with a singular source are discussed in this paper. 
		The stationary discontinuity induced by the singular source and its coupling with the convection of fluid presents challenges to numerical methods.
		We show that the splitting scheme is not well-balanced and leads to incorrect results;
		in addition, some popular well-balanced schemes also give incorrect solutions in extreme cases due to the singularity of source.
		To fix such difficulties, we propose a solution-structure based approximate Riemann solver, in which the structure of Riemann solution is first predicted and then its corresponding approximate solver is given.
		The proposed solver can be applied to the calculation of numerical fluxes in a general finite volume method, which can lead to a new well-balanced scheme.
		Numerical tests show that the discontinuous Galerkin method based on the present approximate Riemann solver has the ability to capture each wave accurately.
	\end{abstract}
\end{abstract}
	
\noindent{\bf Keywords: } hyperbolic conservation law, Euler equation, Riemann problem, singular source, Riemann solver, well-balanced scheme

\noindent{\bf AMS subject classifications: }35L867,35Q31,65M60,76N30

\section{Introduction}\label{introduction}
The governing equations for many mechanical problems involving fluids are hyperbolic conservation laws with singular sources.
For example, surface tension \cite{fechter2018approximate,houim2013ghost}, evaporation\cite{das2020sharp,das2020simulation,das2021sharp}, phase transition\cite{hitz2020comparison,lee2017sharp,lee2017direct,long2021fully,paula2019analysis,suh2008level} and condensation\cite{cheng2010condensation} at the interface of multiphase flows are formulated by singular sources.
Geometric constraints on irregular domains can also lead to the generation of singular sources in simplified governing equations, such as the shallow water equations with discontinuous river beds\cite{george2008augmented,lefloch2011godunov,li2021well} and the governing equations of nozzle flow with discontinuous cross-section\cite{coquel2014robust,thanh2009riemann}. 
In addition, the damping term in the fluid under consideration behaves as a singular source when it is smaller than the grid size\cite{GOSSE2008A}.

In this paper we consider the Euler equations with a singular source:
\begin{equation}\label{equ: governing equation}
	\partial _tU+\partial _xF(U)=\delta(x)S,
\end{equation}
where
\[
U=
\begin{pmatrix}
	\rho\\\rho u\\E
\end{pmatrix},
F=
\begin{pmatrix}
	\rho u\\\rho u^2+p\\(E+p)u
\end{pmatrix},
\delta(x)=\begin{cases}
	0,\ x\ne 0\\
	\infty,\ x=0
\end{cases}
\]
Here, $\rho$, $p$ and $E$ denote the density, pressure and total energy, respectively. $u$ is the velocity. $S$ is the vector of souce term distributed only at the origin. 
Let $U_-=U(0-,t)$ and $U_+=U(0+,t)$ denote the states to the left-hand and right-hand sides of the origin, respectively, and the source leads to $U_-\ne U_+$. 
The source term in this paper is a generalisation of the source term in \cite{cheng2010condensation}:
\begin{equation}\label{equation: value of source}
	S(U_-,U_+)=\begin{cases}
		S(F(U_-))=diag(k_1,k_2,k_3)F(U_-),\ \text{if}\ u_->0,u_+>0\\
		S(F(U_+))=diag(k_1,k_2,k_3)F(U_+),\ \text{if}\ u_-<0,u_+<0\\
		0,\ else\\
	\end{cases}
\end{equation}
where $diag(k_1,k_2,k_3)$ is the diagonal matrix with constant diagonal element $k_1>-1$, $k_2>-1$ and $k_3>-1$. 
This paper deals only with the ideal gas, whose equation of state is 
\begin{equation*}
	R T=p/\rho,\	p=(\gamma-1)\rho e,\ \gamma>1, 
\end{equation*}
where $\gamma$ is the ratio of specific heats, $e$ is the internal energy, $T$ is the temperature and $R$ is the universal gas constant. 

The aim of this paper is to develop a numerical method for hyperbolic conservation laws with singular sources, which is capable of giving numerical solutions in good agreement with the exact solution.

The most commonly used schemes are basesd on splitting techniques, which solves iteratively the partial differential equation associated to the convection and the ordinary differential equation associated to the source (see \cite{1994Upwind,toro2013riemann}). 
Unfortunately, the numerical solution obtained by such a simple method is inaccurate, even if the source term is discretized by a upwind manner.
Small errors in the approximate solution near the singular source may become uncontrolled under the splitting scheme, and this error does not decrease as the grid size decreases, see \cite{audusse2004fast,bernetti2008exact,kroner2008minimum,kroner2005numerical,leroux1999riemann}.
One reason for such error is that the splitting scheme is not well-balanced (also known as C-property, see \cite{abgrall2018high,greenberg1997analysis,li2021well,1994Upwind}), which refers to the ability of the numerical scheme to maintain precisely the equilibrium state near singular source or some steady solutions.
Many scholars have been working on the development of well-balanced schemes, such as \cite{abgrall2018high,audusse2004fast,bernetti2008exact,botchorishvili2003finite,botchorishvili2003equilibrium,castro2007well,castro2008well,leroux1999riemann,toro2013riemann,kroner2005numerical,thanh2021dimensional,gomez2021high,karlsen2009well,leveque1990study,leveque1998balancing,bale2003wave,papalexandris1997unsplit,jin2001steady,jin2005two,2006High}.
Among those schemes, we focus on the well-balanced scheme proposed by Kr{\"o}ner and Thanh in \cite{kroner2005numerical} and referred to in this paper as the K-T scheme to distinguish it from the well-balanced preperty.
One advantage of K-T scheme is that it can be directly generalised to arbitrary equations as long as the solution to the discontinuity induced by singular sources can be obtained explicitly.
We apply K-T scheme to equation (\ref{equ: governing equation}) and obtained good results in the non-extreme case. In some extreme cases, referred to as choked solutions in this paper, the results of K-T scheme are not satisfactory. 
Similar errors of K-T scheme have been reported in \cite{thanh2013numerical}.

The Godunov scheme, which constructs the numerical flux by the exact Riemann solution, is well-balanced for hyperbolic conservation laws with singular sources(see \cite{cuong2015godunov,greenberg1996well,lefloch2011godunov,leroux1999riemann}).
However, the complete Riemann solution of (\ref{equ: governing equation}) is generally not available.
In this paper we apply the approximate Riemann solver to construct the numerical flux at singular sources.
In order for the numerical scheme to be well-balanced, the approximate Riemann solver must give an exact solution when the initial value is a single stationary wave.
Based on our previous analysis of the Riemann problem in \cite{yuTARP}, we propose a structure-based approximate Riemann solver.
The proposed solver not only ensures the well-balanced property of the numerical scheme, but also is effective in extreme cases due to the coupled treatment of sources and convection.
In this paper the numerical scheme of combining the proposed Riemann solver with a discontinuous Galerkin method gives satisfactory numerical solutions in both extreme and non-extreme cases.

This paper is organized as follows. 
Section\ref{section: Preliminaries} is about preliminaries to the Riemann solution.
InSection \ref{Numerical scheme}, we will present splitting methods and unsplitting methods, and their performance for the well-balanced property.
In Section\ref{Numerical flux for singular source} we will be concerned with the numerical flux for singular sources. We will first discuss the numerical flux of the K-T scheme and then present an approximate Riemann solver and the numerical flux based on the solver.
The results of the splitting scheme, the K-T scheme and the present solver-based scheme in conjunction with the discontinuous Galerkin method for solving a steady solution and seven Riemann solutions of different structures, respectively, will be put in Section\ref{Numerical tests}. 
Finally, some conclusions will be given in Section\ref{Conclusions}.

\section{Preliminaries}\label{section: Preliminaries}
In this section we list some conclusions on the following Riemann problems:
\begin{equation}\label{equ: Riemann problem 2}
	\begin{cases}
		\partial _tU+\partial _xF(U)=\delta(x)S \\
		U(x,0)=\begin{cases}
			U_L,\ x<0 \\
			U_R,\ x>0
		\end{cases}
	\end{cases}
\end{equation}

First of all, we give the definition of the stationary wave induced by the singular source in Section\ref{subsection: Stationary wave}.
Then, we give a brief overview of all the elementary waves of the Riemann problem of the classical Euler equations in Section\ref{subsection: other elementary waves}.
Finally, we give all possible structures of Riemann solution in Section\ref{subsection: Riemann problem}.
For the related proofs, readers are referred to \cite{yuTARP}.

\subsection{Stationary wave}\label{subsection: Stationary wave}

The stationary wave is a discontinuity in the Riemann solution that lies constantly on the t-axis.
We denote the left-hand and right-hand states of a stationary wave as $U_-$ and $U_+$, respectively.
Let $u_-$ and $u_+$ denote the velocity of $U_-$ and $U_+$ respectively, then $u_-$ and $u_+$ have the same sign.
In this section we assume $u_->0$ and $u_+>0$.
$U_-$ and $U_+$ satisfy the jump relation:
\begin{equation}\label{equ: jump relation 1}
	F(U_-)+S(F(U_-))=F(U_+),
\end{equation}
therefore the stationary wave is a local steady solution. $U_-$ and $U_+$ are called equilibrium states.

For a given $U_-$, there are two branches of $U_+$ that satisfy (\ref{equ: jump relation 1}) with $k_1,k_2,k_3$ as parameters, which are referred to as subsonic branch and supersonic branch, respectively.
In order to select physical equilibrium states, we employ the following criterion proposed in \cite{yuTARP}:
\begin{Monotonicity Criterion}\label{criterion1}
	Given coefficients $k_1$, $k_2$ and $k_3$, the left-hand state $U_-$ and right-hand state $U_+$ of a stationary wave satisfies $\lambda_k(U_-)\cdot \lambda_k(U_+)\geq 0,\forall 1\leq k\leq 3$.
\end{Monotonicity Criterion}

Let $M_-$ and $M_+$ denote the Mach numbers of $U_-$ and $U_+$, respectively. The admissible sets of $M_-$ and $M_+$ are denoted as $\Gamma_{-}$ and $\Gamma_{+}$, respectively.
Within the restriction of Criterion\ref{criterion1}, $\Gamma_{-}$ and $\Gamma_{+}$ are shown in Table\ref{table: left and right Mach numbers}.

\begin{table}[htbp]
	\centering
	\caption{The admissible sets of $M_-$ and $M_+$.}
	\label{table: left and right Mach numbers}
	\begin{tabular}{ccccc}
		\toprule
		\specialrule{0em}{3pt}{2pt}
		\multirow{2}*{$k$}&\multicolumn{2}{c}{subsonic branch}&\multicolumn{2}{c}{supersonic branch}\\
		\specialrule{0em}{1pt}{2pt}
		&$\Gamma_{-}$& $\Gamma_{+}$&$\Gamma_{-}$&$\Gamma_{+}$\\
		\specialrule{0em}{3pt}{2pt}	
		\midrule
		\specialrule{0em}{1pt}{2pt}
		$(0,+\infty)$&$(0,M_1^*]$&$(0,1]$&$\left[M^*_2,+\infty \right) $&$\left[1,M_3^{*} \right) $\\
		\specialrule{0em}{1pt}{2pt}
		${0}$&$(0,1]$&$(0,1]$&$\left[1,+\infty \right) $&$\left[1,+\infty \right) $\\
		\specialrule{0em}{1pt}{2pt}
		$(-\infty,0)$&$(0,1]$&$(0,M_1^{**}]$&$\left[1,M_3^{**} \right)$&$\left[M_2^{**},+\infty \right)$\\
		\specialrule{0em}{3pt}{2pt}
		\bottomrule
	\end{tabular}
\end{table}
The definitions of the symbols in Table\ref{table: left and right Mach numbers} are given below.
\begin{equation}
	k=(1+k_1)(1+k_3)/(1+k_2)^2-1. 
\end{equation}
\begin{equation}\label{equ: critical Mach number definition1}
	\begin{aligned}
		&M_1^*=\begin{cases}
			\sqrt{\frac{k\gamma+k+1- (\gamma+1)\sqrt{k(k+1)}}{k+1-k\gamma^2}},\ &\text{if}\ k \ne \frac{1}{\gamma^2-1},\ k>0\\
			\sqrt{\frac{\gamma-1}{2\gamma}},\ &\text{if}\ k = \frac{1}{\gamma^2-1}
		\end{cases} \\
		&M_2^*=\begin{cases}
			\sqrt{\frac{k\gamma+k+1+ (\gamma+1)\sqrt{k(k+1)}}{k+1-k\gamma^2}},\ &\text{if}\ 0<k< \frac{1}{\gamma^2-1}\\
			+\infty, \ &\text{if}\ k\geq \frac{1}{\gamma^2-1}
		\end{cases}\\
		&M_3^*=\begin{cases}
			\sqrt{\frac{\gamma+\sqrt{1-k(\gamma^2-1)}}{\gamma-\gamma\sqrt{1-k(\gamma^2-1)}}}, \ &\text{if}\ 0<k< \frac{1}{\gamma^2-1}\\
			+\infty, \ &\text{if}\ k\geq \frac{1}{\gamma^2-1}
		\end{cases}
	\end{aligned}
\end{equation}
\begin{equation}\label{equ: critical Mach number definition2}
	\begin{aligned}
		&M_1^{**}=\sqrt{\frac{1-\sqrt{-k}}{1+\gamma\sqrt{-k}}},\ \text{if}\ k<0\\
		&M_2^{**}=\begin{cases}
			\sqrt{\frac{1+\sqrt{-k}}{1-\gamma\sqrt{-k}}},\ &\text{if}\  -\frac{1}{\gamma^2}<k<0\\
			+\infty, \ &\text{if}\ -1<k\leq -\frac{1}{\gamma^2}
		\end{cases}\\
		&M_3^{**}=\begin{cases}
			\sqrt{\frac{\gamma\sqrt{1+k}+\sqrt{1+k\gamma^2}}{\gamma\sqrt{1+k}-\gamma\sqrt{1+k\gamma^2}}}, \ &\text{if}\ -\frac{1}{\gamma^2}<k<0\\
			1, \ &\text{if}\ -1<k\leq -\frac{1}{\gamma^2}
		\end{cases}
	\end{aligned}
\end{equation}

Given a left-hand state $U_-$ of the stationary wave, the wave curve $\mathcal{D}(U_-,k_1,k_2,k_3)$ denotes all right-hand states $U_+$ that can be connected to $U_-$ by the stationary wave with coefficients $k_1,k_2,k_3$.
Given a right-hand state $U_+$ of the stationary wave, the wave curve $\mathcal{D}^B(U_+,k_1,k_2,k_3)$ denotes all left-hand states $U_-$ that can be connected to $U_+$ by the stationary wave with coefficients $k_1,k_2,k_3$.

Under Criterion\ref{criterion1}, the wave curve of the stationary wave has following properties without proof. For details, readers are referred to \cite{yuTARP}.

\begin{equation}\label{stationary wave curves1}
	\begin{aligned}
		&\mathcal{D}(U_-,k_1,k_2,k_3):\\
		&\begin{aligned}		
			&\rho_+=\rho_-\left(\frac{M_-}{M_+} \right) ^2\frac{\gamma M_+^2+1}{\gamma M_-^2+1}\frac{(1+k_1)^2}{1+k_2},\\
			&u_+=u_-\left(\frac{M_+}{M_-} \right) ^2\frac{\gamma M_-^2+1}{\gamma M_+^2+1}\frac{1+k_2}{1+k_1},\\
			&p_+=p_-\frac{\gamma M_-^2+1}{\gamma M_+^2+1}(1+k_2),
		\end{aligned}\left\lbrace 
		\begin{aligned}
			&\text{if}\ k>0,\ M_+=\begin{cases}
				\sqrt{\frac{1-I}{1+\gamma I}},\ 0<M_-\leq M_1^{*}\\
				\sqrt{\frac{1+I}{1-\gamma I}},\ M_2^{*}\leq M_-<+\infty
			\end{cases}\\
			&\text{if}\ k<0,\ M_+=\begin{cases}
				\sqrt{\frac{1-I}{1+\gamma I}},\ 0<M_-\leq 1\\
				\sqrt{\frac{1+I}{1-\gamma I}},\ 1\leq M_-< M_3^{**}\\
			\end{cases}\\
			&\text{if}\ k=0,\ M_+=\begin{cases}
				\sqrt{\frac{1-I}{1+\gamma I}},\ 0<M_-\leq 1\\
				\sqrt{\frac{1+I}{1-\gamma I}},\ 1<M_-<+\infty\\
			\end{cases}\\	
		\end{aligned}\right. 
	\end{aligned}
\end{equation}
Here, 
\[
I=\frac{\sqrt{(\gamma M_-^2+1)^2-(\gamma+1)M_-^2[(\gamma-1)M_-^2+2](1+k)}}{\gamma M_-^2+1},
\]

\begin{equation}\label{stationary wave curves2}
	\begin{aligned}
		&\mathcal{D}^B(U_+,k_1,k_2,k_3):\\		
		&\begin{aligned}		
			&\rho_-=\rho_+\left(\frac{M_+}{M_-} \right) ^2\frac{\gamma M_-^2+1}{\gamma M_+^2+1}\frac{1+k_2}{(1+k_1)^2},\\
			&u_-=u_+\left(\frac{M_-}{M_+} \right) ^2\frac{\gamma M_+^2+1}{\gamma M_-^2+1}\frac{1+k_1}{1+k_2},\\
			&p_-=p_+\frac{\gamma M_+^2+1}{\gamma M_-^2+1}\frac{1}{1+k_2},
		\end{aligned}\left\lbrace 
		\begin{aligned}
			&\text{if}\ k>0,\ M_-=\begin{cases}
				\sqrt{\frac{1-I^B}{1+\gamma I^B}},\ 0<M_+\leq 1\\
				\sqrt{\frac{1+I^B}{1-\gamma I^B}},\ 1\leq M_+< M^{*}_3\\
			\end{cases}\\
			&\text{if}\ k<0,\ M_-=\begin{cases}
				\sqrt{\frac{1-I^B}{1+\gamma I^B}},\ 0<M_+\leq M^{**}_1\\
				\sqrt{\frac{1+I^B}{1-\gamma I^B}},\ M^{**}_2\leq M_+<+\infty
			\end{cases}\\
			&\text{if}\ k=0,\ M_-=\begin{cases}
				\sqrt{\frac{1-I^B}{1+\gamma I^B}},\ 0<M_-\leq 1\\
				\sqrt{\frac{1+I^B}{1-\gamma I^B}},\ 1<M_-<+\infty\\
			\end{cases}\\	
		\end{aligned}\right. 
	\end{aligned}
\end{equation}
Here, 
\[
I^B=\frac{\sqrt{(\gamma M_+^2+1)^2-(\gamma+1)M_+^2[(\gamma-1)M_+^2+2]/(1+k)}}{\gamma M_+^2+1},
\]

\begin{Property}\label{property: existence and uniqueness}
	For the stationary wave with given left-hand state $U_-$ and coefficients $k_1,k_2,k_3$,
	\begin{itemize}
		\item[(i)] $D(U_-,k_1,k_2,k_3)$ is a singleton set iff
		\begin{equation*}
			k>0\ \&\ M_1^{*}<M_-<M_2^{*}\quad \text{or}\quad 
			k<0\ \&\ M_-\geq M_3^{**};
		\end{equation*}
		\item [(ii)] $D(U_-,k_1,k_2,k_3)$ is a set of two elements iff $k<0\ \&\ M_-=1$;
		\item[(iii)] $D(U_-,k_1,k_2,k_3)=\emptyset$ in other cases.
	\end{itemize}
\end{Property}

An extreme case of the stationary wave is known as the choke, which is defined as
\begin{Definition}
	A stationary wave is said to be choked if its left-hand state $U_-$ and right-hand state $U_+$ satisfy $\exists 1\leq k\leq m, s.t.\ \lambda_k(U_-)\cdot \lambda_k(U_+)=0$.
\end{Definition}
A Riemann solution is said to be choked if it contains a choked stationary wave. 

By (\ref{stationary wave curves1}) and (\ref{stationary wave curves2}), the strength of the stationary wave is determined by $M_-$.
For the convenience of following discussion, we use some notations to represent the solution of stationary wave:
\begin{equation*}
	M_+=g_M^{\mathcal{D}}(M_-),\ 
	\frac{\rho_+}{\rho_-}=g_d^{\mathcal{D}}(M_-),\ 
	\frac{u_+}{u_-}=g_u^{\mathcal{D}}(M_-),\ \frac{p_+}{p_-}=g_p^{\mathcal{D}}(M_-).
\end{equation*}
When $k<0\& M_-=1$, the functions $g_M^{\mathcal{D}}$, $g_d^{\mathcal{D}}$, $g_u^{\mathcal{D}}$, $g_p^{\mathcal{D}}$ are all double-valued, where one is supersonic and the other is subsonic. We need to specify which branch to use in this case.

\subsection{Shock waves, rarefaction waves and contact discontinuities}\label{subsection: other elementary waves}

There are three elementary waves in the Riemann problem for the classical Euler equations, namely Lax shock waves, rarefaction waves and contact discontinuities.
We define some functions to represent the relations between the left-hand and right-hand states of these three elementary waves, and for a detailed introduction readers are referred to \cite{toro2013riemann,chang1989riemann}.

The characteristic domains of $\lambda_1$ and $\lambda_3$ are both genuinely nonlinear. The elementary wave associated with these two characteristic domains is either a shock wave or a rarefaction wave. We denote

\begin{equation}
	\begin{aligned}
		&f_d^{1-\mathcal{S}}(U_0,p)\mathop{=}\limits^{def}\rho_0\frac{(\gamma-1)p_0+(\gamma+1)p}{(\gamma-1)p+(\gamma+1)p_0},\ p\geq p_0\\
		&f_u^{1-\mathcal{S}}(U_0,p)\mathop{=}\limits^{def}u_0-\sqrt{2}(p-p_0)/\sqrt{\rho_0((\gamma+1)p+(\gamma-1)p_0)},\ p\geq p_0\\
		&f_M^{1-\mathcal{S}}(U_0,p)\mathop{=}\limits^{def}f_u^{1-\mathcal{S}}(U_0,p)\sqrt{f_d^{1-\mathcal{S}}(U_0,p)}/\sqrt{\gamma p},\ p\geq p_0\\
		&f_d^{3-\mathcal{S}}(U_0,p)\mathop{=}\limits^{def}\rho_0\frac{(\gamma-1)p_0+(\gamma+1)p}{(\gamma-1)p+(\gamma+1)p_0},\ p\geq p_0\\
		&f_u^{3-\mathcal{S}}(U_0,p)\mathop{=}\limits^{def}u_0+\sqrt{2}(p-p_0)/\sqrt{\rho_0((\gamma+1)p+(\gamma-1)p_0)},\ p\geq p_0\\		
		&f_M^{3-\mathcal{S}}(U_0,p)\mathop{=}\limits^{def}f_u^{3-\mathcal{S}}(U_0,p)\sqrt{f_d^{3-\mathcal{S}}(U_0,p)}/\sqrt{\gamma p},\ p\geq p_0
	\end{aligned}
\end{equation}
Given the left-hand state $U_0$ of a $\lambda_1-$shock wave, the right-hand state with the pressure $p$ satisfies
\begin{equation*}
	\rho=f_d^{1-\mathcal{S}}(U_0,p),u=f_u^{1-\mathcal{S}}(U_0,p),M=f_M^{1-\mathcal{S}}(U_0,p).
\end{equation*}
Given the right-hand state $U_0$ of a $\lambda_3-$shock wave, the left-hand state with the pressure $p$ satisfies
\begin{equation*}
	\rho=f_d^{3-\mathcal{S}}(U_0,p),u=f_u^{3-\mathcal{S}}(U_0,p),M=f_M^{3-\mathcal{S}}(U_0,p).
\end{equation*}

Let
\begin{equation}
	\begin{aligned}
		&f_d^{1-\mathcal{R}}(U_0,p)\mathop{=}\limits^{def}\rho_0\left( p/p_0\right) ^{\frac{1}{\gamma}},\ p\leq p_0\\
		&f_u^{1-\mathcal{R}}(U_0,p)\mathop{=}\limits^{def}u_0-2a_0/(\gamma-1)\left[ \left( p/p_0\right) ^{\frac{\gamma-1}{2\gamma}}-1\right] ,\ p\leq p_0\\
		&f_d^{3-\mathcal{R}}(U_0,p)\mathop{=}\limits^{def}\rho_0\left( p/p_0\right) ^{\frac{1}{\gamma}},\ p\leq p_0\\
		&f_u^{3-\mathcal{R}}(U_0,p)\mathop{=}\limits^{def}u_0+2a_0/(\gamma-1)\left[ \left( p/p_0\right) ^{\frac{\gamma-1}{2\gamma}}-1\right] ,\ p\leq p_0
	\end{aligned}
\end{equation}
Given the left-hand state $U_0$ of a $\lambda_1-$rarefaction wave, the right-hand state with the pressure $p$ satisfies
\begin{equation*}
	\rho=f_d^{1-\mathcal{S}}(U_0,p),u=f_u^{1-\mathcal{S}}(U_0,p).
\end{equation*}
Given the right-hand state $U_0$ of a $\lambda_3-$rarefaction wave, the left-hand state with the pressure $p$ satisfies
\begin{equation*}
	\rho=f_d^{3-\mathcal{S}}(U_0,p),u=f_u^{3-\mathcal{S}}(U_0,p).
\end{equation*}

We denote
\begin{equation}
	\begin{aligned}
		&g_d^{1-\mathcal{R}}(M_0,M)\mathop{=}\limits^{def}g_d^{3-\mathcal{R}}(M_0,M)\mathop{=}\limits^{def}\left[\frac{(\gamma-1)M_0+2}{(\gamma-1)M+2} \right] ^{\frac{2}{\gamma-1}},\ M\geq M_0.\\
		&g_u^{1-\mathcal{R}}(M_0,M)\mathop{=}\limits^{def}g_u^{3-\mathcal{R}}(M_0,M)\mathop{=}\limits^{def}\frac{M}{M_0}\left[\frac{(\gamma-1)M_0+2}{(\gamma-1)M+2} \right] ,\ M\geq M_0.\\
		&g_p^{1-\mathcal{R}}(M_0,M)\mathop{=}\limits^{def}g_p^{3-\mathcal{R}}(M_0,M)\mathop{=}\limits^{def}\left[\frac{(\gamma-1)M_0+2}{(\gamma-1)M+2} \right] ^{\frac{2\gamma}{\gamma-1}} ,\ M\geq M_0.\\
	\end{aligned}
\end{equation}
Given the left-hand state $U_0$ of a $\lambda_1-$rarefaction wave, the right-hand state with the Mach number $M$ satisfies
\begin{equation*}
	\rho=\rho_0g_d^{1-\mathcal{R}}(M_0,M),u=u_0g_u^{1-\mathcal{R}}(M_0,M),p=p_0g_p^{1-\mathcal{R}}(M_0,M).
\end{equation*}
Given the right-hand state $U_0$ of a $\lambda_3-$rarefaction wave, the left-hand state with the Mach number $M$ satisfies
\begin{equation*}
	\rho=\rho_0g_d^{3-\mathcal{R}}(M_0,M),u=u_0g_u^{3-\mathcal{R}}(M_0,M),p=p_0g_p^{3-\mathcal{R}}(M_0,M).
\end{equation*}

The characteristic domain of $\lambda_2=u$ is linearly degenerate. The elementary wave associated with the $\lambda_2$ characteristic domain is the contact discontinuity.

In summary, we have
\begin{equation*}
	\begin{aligned}
		f_d^{1-\mathcal{W}}(U_0,p)=\begin{cases}
			f_d^{1-\mathcal{S}}(U_0,p),\ p\geq p_0\\
			f_d^{1-\mathcal{R}}(U_0,p),\ p<p_0
		\end{cases}
		f_u^{1-\mathcal{W}}(U_0,p)=\begin{cases}
			f_u^{1-\mathcal{S}}(U_0,p),\ p\geq p_0\\
			f_u^{1-\mathcal{R}}(U_0,p),\ p<p_0
		\end{cases}\\
		f_d^{3-\mathcal{W}}(U_0,p)=\begin{cases}
			f_d^{3-\mathcal{S}}(U_0,p),\ p\leq p_0\\
			f_d^{3-\mathcal{R}}(U_0,p),\ p>p_0
		\end{cases}
		f_u^{3-\mathcal{W}}(U_0,p)=\begin{cases}
			f_u^{3-\mathcal{S}}(U_0,p),\ p\leq p_0\\
			f_u^{3-\mathcal{R}}(U_0,p),\ p>p_0
		\end{cases}
	\end{aligned}
\end{equation*}

\subsection{Riemann problem}\label{subsection: Riemann problem}

\begin{figure}[htbp]\label{figure: general structure}
	\centering
	\subfigure{
		\includegraphics[width=0.4\linewidth]{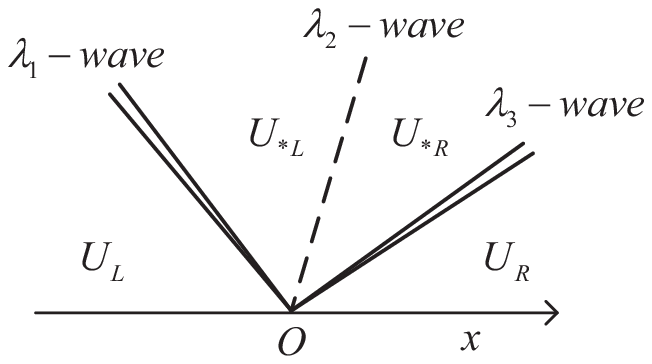}
	}
	\quad
	\subfigure{
		\includegraphics[width=0.4\linewidth]{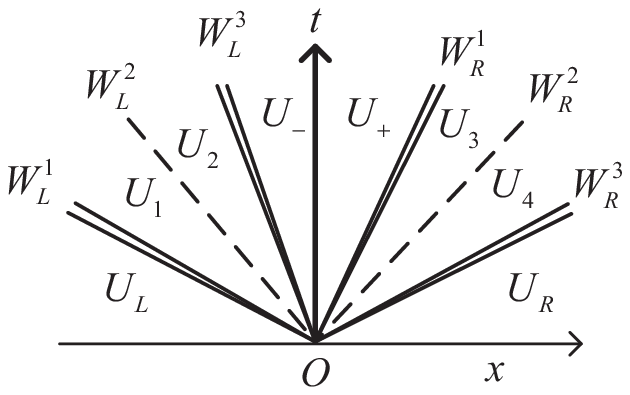}
	}
	\caption{The general structure of Riemann solutions. Left: associate Riemann problem (\ref{equ: Riemann problem 2}). Right: Riemann problem (\ref{equ: associate Riemann problem}).}
\end{figure}

This section presents the structure of the solution to the Riemann problem (\ref{equ: Riemann problem 2}).
We first consider its associate Riemann problem, that is, the following classical Riemann problem for the Euler equations:
\begin{equation}\label{equ: associate Riemann problem}
	\left\lbrace 
	\begin{aligned}
		&\frac{\partial U}{\partial t}+\frac{\partial F}{\partial x}=0,\\
		&U(x,0)=\begin{cases}
			U_L,\ x<0\\
			U_R,\ x>0
		\end{cases}            
	\end{aligned}
	\right. 
\end{equation}\\
The solution of the associative Riemann problem (\ref{equ: associate Riemann problem}) is self-similar and contains four constant regions, which are divided by three classical elementary waves, as shown in the left figure of Figure\ref{figure: general structure}. Although the governing equation (\ref{equ: governing equation}) contains source terms, it is clearly still self-similar, we therefore assume that the solution to the Riemann problem (\ref{equ: Riemann problem 2}) is self-similar.

The self-similar solution of Riemann problem (\ref{equ: Riemann problem 2}) and associate Riemann problem (\ref{equ: associate Riemann problem}) with the initial value condition $U_L$ and $U_R$ are denoted as $U^R\left( \frac{x}{t},U_L,U_R\right),\ U^{A-R}\left( \frac{x}{t},U_L,U_R\right)$, respectively. We have following conclusions as shown in \cite{yuTARP}.

\begin{Proposition}[Double CRP framework]\label{proposition: double CRPs framework}
	The self-similar solution of the Riemann problem consists of seven discontinuities at most. They are a stationary wave at $x=0$, two genuinely nonlinear waves and a contact discontinuity left to $x=0$, two genuinely nonlinear waves and a contact discontinuity right to $x=0$ respectively.
\end{Proposition}
According to Proposition\ref{proposition: double CRPs framework}, the general structure of Riemann solutions is shown in the right figure of Figure\ref{figure: general structure}. The elementary waves on the left and right sides are denoted as $W_L^1$, $W_L^2$, $W_L^3$ and $W_R^1$, $W_R^2$, $W_R^3$, respectively. 
$W_L^1$, $W_L^3$, $W_R^1$ and $W_R^3$ can be shock waves or rarefaction waves, and both $W_L^2$ and $W_R^2$ are contact discontinuities. The t-axis is a stationary wave and $U_-\ne U_+$.
We follow the notation of \cite{lefloch2003riemann,thanh2009riemann} to express the structure of Riemann solutions. For example, $\mathcal{S}(U_L,U_1)$ and $\mathcal{R}(U_L,U_1)$ mean that two states $U_L$ and $U_1$ are connected by a shock and a rarefaction wave, respectively. The symbol $\mathcal{W}(U_L,U_1)$ means that  $U_L$ and $U_1$ are connected by a shock wave or a rarefaction wave.
$\mathcal{C}(U_1,U_2)$ means that $U_1$ and $U_2$ are connected by a contact discontinuity, and $\mathcal{D}(U_-,U_+)$ means that $U_-$ and $U_+$ are connected by a stationary wave. The connection between multiple elementary waves and constant regions is represented by the symbol "$\oplus $". For example, $\mathcal{S}(U_L,U_1)\oplus  \mathcal{C}(U_1,U_2)$ means that $U_L$ and $U_1$ are connected by a shock and a shock wave, followed by a contact discontinuity with the left-hand state $U_1$ and the right-hand state $U_2$.

\begin{Theorem}\label{theorem: structures}
	There are seven structures of the Riemann solution:
	\begin{itemize}
		\item [(i)]non-choked structures:
		\begin{itemize}
			\item[] Type1: $\mathcal{W}(U_L,U_-)\oplus \mathcal{D}(U_-,U_+)\oplus \mathcal{C}(U_+,U_4)\oplus \mathcal{W}(U_4,U_R)$.
			\item[] Type2: $\mathcal{D}(U_L,U_+)\oplus \mathcal{W}(U_+,U_3)\oplus \mathcal{C}(U_3,U_4)\oplus \mathcal{W}(U_4,U_R)$.
		\end{itemize}
		\item [(ii)]choked structures:
		\begin{itemize}
			\item[] Type3: $\mathcal{W}(U_L,U_-)\oplus \mathcal{D}(U_-,U_+)\oplus \mathcal{R}(U_+,U_3)\oplus \mathcal{C}(U_3,U_4)\oplus \mathcal{W}(U_4,U_R)$ with $M_+=1$.
			\item[] Type4: $\mathcal{D}(U_L,U_+)\oplus \mathcal{R}(U_+,U_3)\oplus \mathcal{C}(U_3,U_4)\oplus \mathcal{W}(U_4,U_R)$ with $M_+=1$.
			\item[] Type5: $\mathcal{R}(U_L,U_-)\oplus \mathcal{D}(U_-,U_+)\oplus \mathcal{W}(U_+,U_3)\oplus \mathcal{C}(U_3,U_4)\oplus \mathcal{W}(U_4,U_R)$ with $M_-=1$.
			\item[] Type6: $\mathcal{R}(U_L,U_-)\oplus \mathcal{D}(U_-,U_+)\oplus \mathcal{C}(U_+,U_4)\oplus \mathcal{W}(U_4,U_R)$ with $M_-=1$.
			\item[] Type7: $\mathcal{R}(U_L,U_-)\oplus \mathcal{D}(U_-,U_+)\oplus \mathcal{R}(U_+,U_3)\oplus \mathcal{C}(U_3,U_4)\oplus \mathcal{W}(U_4,U_R)$ with $M_-=M_+=1$.
		\end{itemize}
	\end{itemize}		
\end{Theorem}

\begin{Theorem}\mbox{}
	\begin{itemize}
		\item[(i)] If $k>0$, then all possible structures are Type1, Type2, Type3 and Type4.
		\item[(ii)] If $k<0$, then all possible structures are Type1, Type2, Type5 and Type6.
		\item[(iii)] If $k=0$, then all possible structures are Type1, Type2 and Type7.	 
	\end{itemize}
\end{Theorem}

\begin{Remark}\label{remark: type4}
	The initial value of the Riemann solution of Type4 satisfies $M_L=M^*_2$. The Riemann solution of Type4 can be regarded as a Riemann solution of Type3 for which the wave $W_L^1$ is a stationary shock, and it can also be regarded as a Riemann solution of Type2 with $M_L=M^*_2$. 
\end{Remark}

\begin{Remark}\label{remark: type6}
	The initial value of the Riemann solution of Type6 satisfies
	\begin{equation}\label{type6 initial value}
		u_Lg_u^{1-\mathcal{R}}(M_L,1)g_u^{\mathcal{D}}(1)-f_u^{3-\mathcal{W}}(U_R,p_Lg_p^{1-\mathcal{R}}(M_L,1)g_p^{\mathcal{D}}(1))=0,
	\end{equation}
	where $g^{\mathcal{D}}_u(1)$ and $g^{\mathcal{D}}_p(1)$ are both supersonic branches. The Riemann solution of Type6 can be regarded as a Riemann solution of Type5 for which the wave $W_R^1$ is a stationary shock, and it can also be regarded as a Riemann solution of Type1 with (\ref{type6 initial value}). 
\end{Remark}

\section{Numerical schemes}\label{Numerical scheme}
In this section we focus on the numerical discretization of the inhomogeneous system (\ref{equ: governing equation}). The numerical computation is carried out for the following Cauchy problem:
\begin{equation}\label{equ: Cauchy problem}
	\begin{cases}
		\partial _tU+\partial _xF(U)=\delta(x)S,\ x\in\left(-\infty,+\infty \right) ,\ t> 0\\
		U(x,0+)=U_0(x).
	\end{cases}
\end{equation}

The computational domain $\left[a,b \right] (a<0,b>0)$ is discretized using $N$ non-overlapping cells with cell interfaces $a=x_{\frac{1}{2}}<x_{\frac{3}{2}}<...<x_{N+\frac{1}{2}}=b$. The center of cell $I_j=\left[ x_{j-\frac{1}{2}},x_{j+\frac{1}{2}} \right] $ is denoted by $x_j=(x_{j-\frac{1}{2}}+x_{j+\frac{1}{2}})/2$. For the rest of this paper, we assume the discretization to be uniform with the mesh size denoted by $h=x_{j+\frac{1}{2}}-x_{j-\frac{1}{2}}$. Given a uniform time step $\Delta t$, setting $t_n=n\Delta t,\ n\in \mathbf{N}$. Let $U_j^n$ represents the cell average of the approximate solution $U^h(x,t)$ at time $t^n$,
\begin{equation}\label{cell average approximation}
	\begin{aligned}
		&U_j^n=\frac{1}{h}\int_{I_j}U^h(x,t^n)dx,\\
		&U_j^0=\frac{1}{h}\int_{I_j}U_0(x)dx.
	\end{aligned}
\end{equation}
The approximate solution $U^h(x,t)$ at time $t^n$ is $U^h(x,t^n)|_{I_j}=U_j^n$.
A general numerical scheme can be written as 
\begin{equation*}
	U_j^{n+1}=\mathcal{H}(U_{j-l}^n,...,U_{j+l}^n).
\end{equation*}

The exact solution of ((\ref{equ: Cauchy problem})) is discontinuous at the origin. 
If the origin lies in the interior of a cell $I_j$, the approximation of a smooth approximate solution on $I_j$ to an discontinuous exact solution inevitably leads to a large numerical dissipation, with the result that the approximation of the value of source is wrong.
A better approach would be to place the origin at the interface of two cells.
Thus we use the following grid:
\begin{equation}\label{grid}
	a=x_{\frac{1}{2}}<x_{\frac{3}{2}}<...<x_{j_0+\frac{1}{2}}=0<...<x_{N+\frac{1}{2}}=b.
\end{equation}

The well-balanced property is related to the numerical approximation of steady solutions. We focus on the discrete methods of singular source term, therefore we need to consider whether the numerical scheme can maintain a stationary wave.
\begin{Definition}\label{definition: well-balanced}
	Suppose the initial value of ((\ref{equ: Cauchy problem})) is
	\begin{equation*}
		U_0(x)=\begin{cases}
			U_-,\ x\leq 0\\
			U_+,\ x>0
		\end{cases}
	\end{equation*}
	A numerical scheme for solving ((\ref{equ: Cauchy problem})) is said to be well-balanced if, given any $U_-$ with $M_-\in \Gamma_-$ and $U_+\in D(U_-,k_1,k_2,k_3)$ such that 
	\begin{equation*}
		U_{j}^{n}=U_{j}^0,\ \forall n\geq 0, \forall 1\leq  j\leq N .
	\end{equation*}
\end{Definition}

Under the grid(\ref{grid}) and approximate solution (\ref{cell average approximation}), the initial value in Definition\ref{definition: well-balanced} is equivalent to
\begin{equation}\label{equation: well-balanced initial value}
	U_j^0=\left\lbrace 
	\begin{aligned}
		&U_-,1\leq j\leq j_0\\
		&U_+,j_0+1\leq j\leq N
	\end{aligned}
	\right. 
\end{equation}

The numerical discretization of source can be divided into two types: splitting methods and unsplitting methods. 
In splitting methods the source and convection are discretized separately, and in unsplitting methods the source and the convection are coupled, such as \cite{toro2013riemann,yang2013discontinuous}.
In the following we give the general form of splitting methods and unsplitting methods for (\ref{equ: Cauchy problem}) and some analysis about the well-balanced property.

\subsection{Splitting method}
The splitting method alternates between solving homogeneous Euler equations in one step, and the ordinary differential equation of source terms in the second step.
\begin{subequations}	
	\begin{align}\label{splitting: PDE}
		&\left. 
		\begin{aligned}
			\text{PDE: }&U_t+F(U)_x=0,\\
			\text{IC: }&U^h(x,t^n),
		\end{aligned}
		\right\rbrace \Longrightarrow \bar{U}^h(x,t^{n+1})\\
		\label{splitting: ODE}
		&\left. 
		\begin{aligned}
			\text{ODE: }&\frac{d}{dt}U=\delta(x)S,\\
			\text{IC: }&\bar{U}^h(x,t^{n+1}),
		\end{aligned}
		\right\rbrace \Longrightarrow U^h(x,t^{n+1})
	\end{align}
\end{subequations}

We denote by $C^{(t)}$ and $S^{(t)}$ the solution operators for (\ref{splitting: PDE}) of convection and (\ref{splitting: ODE}) of source for a time $t$, respectively, then the splitting method can be expressed as
\begin{equation}
	U^{n+1}=S^{(\Delta t)}C^{(\Delta t)}(U^n).
\end{equation}

\begin{Theorem}\label{thoerem: splitting}
	The splitting method is not well-balanced for any given $\Delta t>0$.
\end{Theorem}
\begin{Proof}
	Let $U_0(x)$ be the initial value given in (\ref{equation: well-balanced initial value}).
	For the splitting method to be well-balanced, it must satisfy
	\begin{equation}\label{splitting well-balanced property}
		S^{(\Delta t)}C^{(\Delta t)}(U_0)=U_0.
	\end{equation}
	In the following we assume that both $C^{(t)}$ and $S^{(t)}$  are exactly solved, then we prove that (\ref{splitting well-balanced property}) is not valid.
	
	First we prove
	\begin{equation}\label{splitting theorem: conclusion}
		\|C^{(\Delta t)}(U_0)-U_0\|_{L^2}\ne 0.
	\end{equation}
	The Cauchy problem for the convection equation with $U_0(x)$ as the initial value is a classical Riemann problem for the Euler equations: $CRP(U_-,U_+)$.
	If (\ref{splitting theorem: conclusion}) was false, then the solution to $CRP(U_-,U_+)$ is a single stationary contact discontinuity or a single stationary shock wave. For the case of contact discontinuity we have $u_-=u_+=0$, which is in contradiction with $U_-\ne U_+$. 
	For the case of shock wave we have $(M_--1)(M_+-1)<0$, which is in contradiction with Criterion\ref{criterion1}.
	
	Second, we show that $S^{(\Delta t)}\left( C^{(\Delta t)}(U_0)\right) $ does not return to the initial equilibrium state. 
	By (\ref{splitting theorem: conclusion}), it follows that there exists $\xi >0$ such that $C^{(\Delta t)}(U_0)(\xi)\ne 0$.
	The ODE at $x=\xi$ is $\frac{d}{dt}U(\xi)=\delta(\xi)S=0$, therefore \[S^{(\Delta t)}\left( C^{(\Delta t)}(U_0)\right)(\xi)=C^{(\Delta t)}(U_0)(\xi)\ne 0.\]
	
	We have thus proved $S^{(\Delta t)}C^{(\Delta t)}(U_0)\ne U_0$.
\end{Proof}

For the numerical discretization of the PDE, we can choose one from a large number of existing consistent and stable schemes. However, Theorem\ref{thoerem: splitting} shows that no choice can make the splitting method well-balanced.
For the numerical discretization of the ODE, we give a upwind $S^{(t)}$ as follows.
\begin{equation}\label{splitting method-ODE}\left\lbrace 
	\begin{aligned}
		&\text{if }u_{j_0}, u_{j_0+1}> 0,\ U_j^{n+1}
		=\begin{cases}
			U_j^n+\frac{\Delta t}{h}S(F( U_{j_0}^n )),\ &j=j_0+1\\
			U_j^{n},\ &j\ne j_0+1
		\end{cases}\\
		&\text{if }u_{j_0}, u_{j_0+1}< 0,\ U_j^{n+1}
		=\begin{cases}
			U_j^n+\frac{\Delta t}{h}S(F( U_{j_0+1}^n )) ,\ 	&j=j_0\\
			U_j^{n},\ &j\ne j_0
		\end{cases}\\
		&\text{else, }U_j^{n+1}=U_j^{n}
	\end{aligned}\right. 
\end{equation}
It can be shown that the above operator satisfies $S^{(\Delta t)}(U^0)=U^0$ for any $U^0$ that satisfies (\ref{equation: well-balanced initial value}).

\subsection{Unsplitting method}
In order to give a general form of unsplitting scheme, we perform some formal derivations of the numerical discretization of the inhomogeneous system (\ref{equ: Cauchy problem}).

Integrating (\ref{equ: Cauchy problem}) over space $I_j$ and time $t^n\leq t\leq t^{n+1}$, we have
\begin{equation}\label{equation: unsplitting scheme}
	\begin{aligned}
		U^{n+1}_j=&U_j^n-\frac{\Delta t}{h}\left[ \int_{t^n}^{t^{n+1}}F(U(x_{j+\frac{1}{2}}^-,t))dt
		-\int_{t^n}^{t^{n+1}}F(U(x_{j-\frac{1}{2}}^+,t))dt
		\right] \\
		&+\Delta t \int_{t^n}^{t^{n+1}}S(U(x_{j_0-\frac{1}{2}}^-,t),U(x_{j_0+\frac{1}{2}}^+,t))dt
		\int_{I_j}\delta(x)dx,
	\end{aligned}
\end{equation}
where $x^{\pm}=\lim_{\epsilon \to 0}x\pm \epsilon$.

We denote the approximation of the integral of $F$ in $[t^n,t^{n+1}]$ in the above equation separately by
\begin{equation}\label{unsplitting: flux other than origin}
	\begin{aligned}
		&\widehat{F}_{j+\frac{1}{2},-}^n=\widehat{F}_-(U(x_{j+\frac{1}{2}}^-,t^n),U(x_{j+\frac{1}{2}}^+,t^n))
		\approx \int_{t^n}^{t^{n+1}}F(U(x_{j+\frac{1}{2}}^-,t))dt\\
		&\widehat{F}_{j+\frac{1}{2},+}^n=\widehat{F}_+(U(x_{j+\frac{1}{2}}^-,t^n),U(x_{j+\frac{1}{2}}^+,t^n))
		\approx \int_{t^n}^{t^{n+1}}F(U(x_{j-\frac{1}{2}}^+,t))dt
	\end{aligned}
\end{equation}

At the cell interface $x_{j+\frac{1}{2}}(j\ne j_0)$, the equations is locally conservative, so
\begin{equation}
	\widehat{F}_{j+\frac{1}{2},-}^n=\widehat{F}_{j+\frac{1}{2},+}^n.
\end{equation}

If $j\ne j_0$ and $j\ne j_0+1$, then we have
\begin{equation*}
	\int_{I_j}\delta(x)dx=0,
\end{equation*}
and (\ref{equation: unsplitting scheme}) is
\begin{equation}
	U^{n+1}_j=U_j^n-\frac{\Delta t}{h}\left[ \widehat{F}_{j+\frac{1}{2},-}^n-\widehat{F}_{j-\frac{1}{2},+}^n \right] ,\ j\ne j_0,j\ne j_0+1.
\end{equation} 

The $S$ at time $t$ is $S(U(x_{j_0-\frac{1}{2}}^-,t),U(x_{j_0+\frac{1}{2}}^+,t))$, which is independent of the spatial coordinate $x$. We abbreviate it $S(t)$.

In cell $I_{j_0}$, the scheme is 
\begin{equation}
	U^{n+1}_{j_0}=U_{j_0}^n-
	\frac{\Delta t}{h}\left[ \int_{t^n}^{t^{n+1}}
	\left( F(U(x_{{j_0}+\frac{1}{2}}^-,t))
	-h\int_{I_{j_0}}\delta(x)dxS(t)
	\right) dt
	-\widehat{F}_{{j_0}-\frac{1}{2}}^n
	\right].
\end{equation}

For simplicity of notation, we continue to use the letter $\widehat{F}$.
\begin{equation}\label{unsplitting: left flux}
	\widehat{F}_{{j_0}+\frac{1}{2},-}^n=\widehat{F}_-(U(x_{j_0+\frac{1}{2}}^-,t^n),U(x_{j_0+\frac{1}{2}}^+,t^n))
	\approx \int_{t^n}^{t^{n+1}}
	\left( F(U(x_{{j_0}+\frac{1}{2}}^-,t))
	-h\int_{I_{j_0}}\delta(x)dxS(t)
	\right) dt
\end{equation}
Then (\ref{equation: unsplitting scheme}) can still be written in a conservation form for cell $I_{j_0}$.
\begin{equation}
	U^{n+1}_{j_0}=U_{j_0}^n-\frac{\Delta t}{h}\left[ \widehat{F}_{{j_0}+\frac{1}{2},-}^n-\widehat{F}_{{j_0}-\frac{1}{2},+}^n \right] 
\end{equation} 

Similarly, for cell $I_{j_0+1}$, we have
\begin{equation}\label{unsplitting: right flux}
	\begin{aligned}
		&U^{n+1}_{j_0+1}=U_{j_0+1}^n-\frac{\Delta t}{h}\left[ \widehat{F}_{{j_0}+3/2,-}^n-\widehat{F}_{{j_0}+\frac{1}{2},+}^n \right]\\
		&\widehat{F}_{{j_0}+\frac{1}{2},+}^n=\widehat{F}_+(U(x_{j_0+\frac{1}{2}}^-,t^n),U(x_{j_0+\frac{1}{2}}^+,t^n)),
		\approx \int_{t^n}^{t^{n+1}}
		\left( F(U(x_{{j_0}+\frac{1}{2}}^+,t))
		+h\int_{I_{j_0+1}}\delta(x)dxS(t)
		\right) dt.
	\end{aligned}
\end{equation}

The unsplitting scheme is therefore in the form of a conservative scheme:
\begin{equation}\label{unsplitting conservative scheme}
	U^{n+1}_j=U_j^n-\frac{\Delta t}{h}\left[ \widehat{F}_{j+\frac{1}{2},-}^n-\widehat{F}_{j-\frac{1}{2},+}^n \right] .
\end{equation} 

The only thing that remains to be done is to give a definition of the numerical flux $\widehat{F}$.
At the cell interface $x_{j+\frac{1}{2}}(j\ne j_0)$ we can use some common numerical fluxes for Euler equations, such as Lax-Friedrichs flux, Roe flux, etc.
At the cell interface $x_{j_0+\frac{1}{2}}$, because of
\begin{equation*}
	\int_{I_{j_0}}\delta(x)dx \text{ and } \int_{I_{j_0+1}}\delta(x)dx,
\end{equation*}
we need to design the numerical fluxes$\widehat{F}_{{j_0}+\frac{1}{2},-}^n$ and $\widehat{F}_{{j_0}+\frac{1}{2},+}^n$ specifically taking into account the effect of source.

Integrating (\ref{equ: Cauchy problem}) at $[a,b]\times [t^n,t^{n+1}]$, we have
\begin{equation}\label{unsplitting: integrating1}
	\int_{a}^{b}U(x,t^{n+1})dx=	\int_{a}^{b}U(x,t^{n})dx
	-\int_{t^n}^{t^{n+1}}(F(U(b,t))-F(U(a,t)))dt
	+\int_{t^n}^{t^{n+1}} S(t)dt.
\end{equation}

Summing (\ref{unsplitting conservative scheme}) over all cells, we have
\begin{equation}\label{unsplitting: integrating2}
	\sum_jU_j^{n+1}=\sum_jU_j^{n}-\frac{\Delta t}{h}
	\left( \widehat{F}^n_{N+\frac{1}{2},-}-\widehat{F}^n_{\frac{1}{2},+} \right) 
	-\frac{\Delta t}{h}\left( \widehat{F}^n_{j_0+\frac{1}{2},-}-\widehat{F}^n_{j_0+\frac{1}{2},+} \right) 
\end{equation}

Comparing (\ref{unsplitting: integrating1}) and (\ref{unsplitting: integrating2}), we have
\begin{equation}\label{equation: flux difference}
	\widehat{F}_{j_0+\frac{1}{2},+}^n-\widehat{F}_{j_0+\frac{1}{2},-}^n\approx 
	\int_{t^n}^{t^{n+1}}S(t)dt,
\end{equation}
The effect of source is therefore approximated by the difference between $\widehat{F}_{j_0+\frac{1}{2},+}^n$ and $\widehat{F}_{j_0+\frac{1}{2},-}^n$.

We now consider the consistency condition for numerical fluxes. 
For the numerical flux used at the interface $x_{j+\frac{1}{2}}(j\ne j_0)$, it should satisfy
\begin{equation}\label{consistency condition 1}
	\text{consistency condition: } \widehat{F}_-(U,U)=\widehat{F}_+(U,U)=F(U),\ \forall U.
\end{equation}
At the cell interface $x_{j_0+\frac{1}{2}}$, the consistency condition above does not take into account the effect of source and is therefore inappropriate.
We define the following consistency condition.
\begin{equation}\label{consistency condition 2}
	\begin{aligned}
		\text{consistency condition: } &\widehat{F}_-(U_-,U_+)=F(U_-),\widehat{F}_+(U_-,U_+)=F(U_+),\\ 
		&\forall U_- \text{ with } M_-\in \Gamma_- \text{ and } U_+\in D(U_-,k_1,k_2,k_3).
	\end{aligned}
\end{equation}
When $S$ degenerates to zero, the consistency condition (\ref{consistency condition 2}) degenerates to the consistency condition (\ref{consistency condition 1}).

\begin{Theorem}\label{theorem: unsplitting well-balancd}
	Under the grid (\ref{grid}), the conservative scheme (\ref{unsplitting conservative scheme}) is well-balanced if for any $n\geq 0$, the numerical fluxes $\widehat{F}^n_{j+\frac{1}{2}\pm}$ satisfy the consistency condition (\ref{consistency condition 1}) for any $j\ne j_0$ and the consistency condition (\ref{consistency condition 2}) for $j=j_0$.
\end{Theorem}
\begin{Proof}
	Suppose that the initial value of the Cauchy problem (\ref{equ: Cauchy problem}) satisfies the conditions in Definition\ref{definition: well-balanced}. 
	At the cell boundary $x_{j+\frac{1}{2}}(j\ne j_0)$, it follows from the monotonicity condition (\ref{consistency condition 1}) that
	\begin{equation*}
		\begin{aligned}
			&\widehat{F}_{j+\frac{1}{2},-}^0=\widehat{F}_{j+\frac{1}{2},+}^0=F(U_-),\text{ if }j<j_0\\
			&\widehat{F}_{j+\frac{1}{2},-}^0=\widehat{F}_{j+\frac{1}{2},+}^0=F(U_+),\text{ if }j>j_0
		\end{aligned}
	\end{equation*}
	At the cell boundary $x_{j_0+\frac{1}{2}}$, it follows from the monotonicity condition (\ref{consistency condition 2}) that
	\begin{equation*}
		\widehat{F}_{j_0+\frac{1}{2},-}^0=F(U_-), \widehat{F}_{j_0+\frac{1}{2},+}^0=F(U_+).
	\end{equation*}
	With the conservative scheme (\ref{unsplitting conservative scheme}), we have
	\begin{equation*}
		U_j^1= U_j^0,\forall 1\leq j\leq N.
	\end{equation*}
	Therefore we have
	\begin{equation*}
		U_{j}^{n}=U_{j}^0,\ \forall n\geq 0, \forall 1\leq  j\leq N .
	\end{equation*}
\end{Proof}

\section{Numerical flux for singular source}\label{Numerical flux for singular source}

As the splitting method is not well-balanced, we are concerned with the unsplitting method.
As long as the numerical fluxes satisfy the consistency condition, the scheme is well-balanced. However, due to the discontinuity of the wave curve of stationary waves and the restriction of Criterion\ref{criterion1}, a well-balanced scheme may still lead to unstable or even incorrect numerical solutions.

In the following we will first present the numerical flux in the K-T sheme to show the difficulties caused by the singular source,
then propose a solution-structure based approximate Riemann solver and apply it to the numerical fluxes to obtain a new class of well-balanced schemes.

\subsection{K-T numerical flux}

In \cite{kroner2005numerical} Kr{\"o}ner and Thanh design a special numerical flux for the singular source induced by the discontinuous bottom in the shallow water equations. We apply that flux to the singular source in (\ref{equ: governing equation}) and call it the K-T numerical flux.
\begin{equation}\label{W-B: flux}
	\begin{aligned}
		&\widehat{F}_{{j_0}+\frac{1}{2},-}^n=
		F\left( U^{A-R}\left(\frac{x}{t}=0, U(x_{{j_0}+\frac{1}{2}}^-,t),D^B(U(x_{{j_0}+\frac{1}{2}}^+,t),k_1,k_2,k_3) \right) \right) ,\\
		&\widehat{F}_{{j_0}+\frac{1}{2},+}^n=F\left( U^{A-R}\left(\frac{x}{t}=0, D(U(x_{{j_0}+\frac{1}{2}}^-,t),k_1,k_2,k_3), U(x_{{j_0}+\frac{1}{2}}^+,t) \right)\right) .
	\end{aligned} 
\end{equation}
The K-T scheme constructs numerical fluxes based on the solutions of two classical Riemann problems. 
We can also use some numerical fluxes as approximate solutions, such as the Lax-Friedrichs numerical flux:
\begin{equation}\label{W-B: L-F flux}
	\begin{aligned}
		&\widehat{F}_{{j_0}+\frac{1}{2},-}^n=
		\widehat{F}^{LF}\left( U(x_{{j_0}+\frac{1}{2}}^-,t),D^B(U(x_{{j_0}+\frac{1}{2}}^+,t),k_1,k_2,k_3) \right) ,\\
		&\widehat{F}_{{j_0}+\frac{1}{2},+}^n=\widehat{F}^{LF}\left( D(U(x_{{j_0}+\frac{1}{2}}^-,t),k_1,k_2,k_3), U(x_{{j_0}+\frac{1}{2}}^+,t) \right).
	\end{aligned} 
\end{equation}

It is easy to check that (\ref{W-B: flux}) satisfies the consistency condition, therefore by Theorem\ref{theorem: unsplitting well-balancd} we have
\begin{Theorem}
	The K-T scheme is well-balanced for any $h$ and given CFL condition.
\end{Theorem}

The K-T scheme applies the wave curve of the stationary wave explicitly to the construction of the numerical fluxes (\ref{W-B: flux}), which not only ensures the well-balanced property of the numerical scheme, but also reduces the numerical dissipation of the numerical solution at the origin.

In some extreme tests of the shallow water equations, as shown in \cite{thanh2013numerical}, the K-T scheme is unable to give satisfactory results.
For some extreme tests of the Cauchy problem(\ref{equ: Cauchy problem}), such as when the left-hand state of the origin is supersonic and the right-hand state is subsonic, the K-T scheme encounters similar difficulties. We believe there are two reasons for such difficulties.

One reason is the unreasonable choice of the branch of the stationary wave curve. 
There are two branches of the stationary wave curve and it is only reasonable to use the same branch in the numerical fluxes $\widehat{F}_{{j_0}+\frac{1}{2},-}^n$ and $\widehat{F}_{{j_0}+\frac{1}{2},+}^n$.
However, the K-T scheme chooses the branch according to the single-sided state, which leads to the possibility that the branches used in fluxes $\widehat{F}_{{j_0}+\frac{1}{2},-}^n$ and $\widehat{F}_{{j_0}+\frac{1}{2},+}^n$ are different.
The discontinuous property of the stationary wave curve causes the error in the numerical fluxes due to unreasonable selection of branches to become uncontrollable.

Another reason is the unresolvability of the numerical fluxes (\ref{W-B: flux}) due to the admissible range of states on either side of the stationary wave.
Under Criterion\ref{criterion1}, the admissible states on either side of the stationary wave are shown in Table\ref{table: left and right Mach numbers}.
When the flow field is in non-equilibrium, the states on either side of the origin may fall within that range. 
In such extreme cases, $D^B(U(x_{{j_0}+\frac{1}{2}}^+,t),k_1,k_2,k_3)$ and $D(U(x_{{j_0}+\frac{1}{2}}^-,t),k_1,k_2,k_3)$ in the numerical fluxes (\ref{W-B: flux}) are unsolvable.
Specifically, we have
\begin{Theorem}\label{theorem: K-T scheme}
	The numerical fluxes of K-T scheme are unavailable if $\exists n\geq 0$ s.t. $M_-^{h,n}\notin \Gamma_-$ or $M_+^{h,n}\notin \Gamma_+$ for given $h>0$, where $M_-^{h,n}$ and $M_+^{h,n}$ are the Mach numbers of $U^h(x_{j_0+\frac{1}{2}}^-,t^n)$ and $U^h(x_{j_0+\frac{1}{2}}^+,t^n)$ respectively, $\Gamma_-$ and $\Gamma_+$ are defined by Table\ref{table: left and right Mach numbers}.
\end{Theorem}

These two errors cause incorrect numerical fluxes and hence lead to large errors in the numerical solutions.
We believe that this is the main reason for the poor performance of the K-T scheme in some tests in Section\ref{Numerical tests}.

\subsection{Solution-structure based approximate Riemann solver and its well-balanced scheme}
In this section we present a new approximate Riemann solver for (\ref{equ: Riemann problem 2}) and the numerical fluxes based on the proposed solver.

\subsubsection{Solution-structure based approximate Riemann solver}
Given the initial values $U_L$ and $U_R$ of the Riemann problem, the present proposed solver gives an approximate solution $\widetilde{U}_-$ for $U^R(\frac{x}{t}=0-,U_L,U_R)$ and an approximate solution $\widetilde{U}_+$ for $U^R(\frac{x}{t}=0+,U_L,U_R)$.

If $u_L\cdot u_R\leq 0$, then the source (\ref{equation: value of source}) is zero, and the governing equations degenerates to the homogeneous Euler equations.
We use the solution of a $CRP$ as the solution to the approximate Riemann solver:
\begin{equation}\label{new scheme: flux 1}
	\widetilde{U}_-=\widetilde{U}_+=U^{A-R}(\frac{x}{t}=0,U_L,U_R).
\end{equation}
In the following we still assume $u_L>0,u_R>0$, since the solution in the case of $u_L<0,u_R<0$ can be obtained by the transformation $x\mapsto -x,\ u\mapsto -u$.
The solver's strategy is to first predict the structure of the solution based on initial values and then solve for the states on either side of the t-axis based on the predicted structure.
Below we present the implementation of these two steps in turn.

\textbf{Step1: prediction of sturctures: }

There are seven structures of Riemann solutions, as shown in Theorem\ref{theorem: structures}.
By Remark\ref{remark: type4} and Remark\ref{remark: type6}, we ignore the two structures of Type4 and Type6 in the approximate Riemann solver.
The alternative predicted structures are the two non-choked structures of Type1 and Type2 as well as the three choked structures of Type3, Type5 and Type7.
The prediction strategy is to first determine whether the Riemann solution is Type2 and Type1 based on $U_L$ and $U_R$ in turn, and if neither is the case, then determine which choked structure the solution is based on the value of $k$.

In the Type2 structure, the initial left-hand flow goes through the singular source point, and we use two conditions for judgement.
One condition is that $U_L$ should be within the admissible set of left-hand states of stationary waves:
\begin{equation}\label{solver: type2 prediction1}
	\left\lbrace 
	\begin{aligned}
		&k>0:\ M_L\geq M_2^*,\\
		&k=0:\ M_L\geq 1,\\
		&k<0:\ 1\leq M_L<M_3^{**},
	\end{aligned}\right. 
\end{equation}
Another condition is that the speed of elementary wave $W_R^1$ is not less than 0, which is equivalent to
\begin{equation}\label{solver: type2 prediction2}
	p\geq p_+ \Rightarrow (\gamma+1)p\leq (2\gamma M_+^2-\gamma+1)p_+,
\end{equation}
where $p_+$ and $M_+$ are the pressure and Mach number of supersonic $U_+=\mathcal{D}(U_L,k_1,k_2,k_3)$,  respectively, and $p$ is the pressure of the classical Riemann solution $U^{A-R}\left( \frac{x}{t},U_+,U_R \right) $ at the intermediate region. If both conditions can be satisfied, then we predict the structure of Riemann solution to be Type2.

Otherwise, we check whether $U_L$ and $U_R$ satisfy the condition of Type1 structure. Associating all the waves in Type1 structure, we see that the pressure of $U_-$ is a root of function $T(p)$, which is defined by
\begin{equation*}
	T(p)\mathop{=}\limits^{def}f_u^{3-\mathcal{W}}(U_R,pg_p^\mathcal{D}(f_M^{1-\mathcal{W}}(U_L,p)))-f_u^{1-\mathcal{W}}(U_L,p)g_u^\mathcal{D}(f_M^{1-\mathcal{W}}(U_L,p)).
\end{equation*}
Using the monotonicity of function $f_M^{1-W}$ with respect to the second argument, we can define $p_1$ and $p_2$ are the pressures that satisfy
\begin{equation*}
	f_M^{1-W}(U_L,p_2)=\begin{cases}
		M_*^1,\ k>0\\
		1,\ k\leq 0
	\end{cases}
\end{equation*}
and $f_M^{1-W}(U_L,p_1)=0$, respectively.
According to Table\ref{table: left and right Mach numbers}, it follows that $0<M_-\leq M_1^*$ for $k>0$, and $0<M_-\leq 1$ for $k\leq 0$, which is equivalent to $p_2\leq p_-<p_1$. 
If $T(p_1)T(p_2)<0$, then we predict the structure of Riemann solution to be Type1.

If none of these conditions can be met, then we predict the structure of Riemann solution to be a choked structure. Precisely, if $k>0$, then we predict the structure to be Type3; if $k=0$, then we predict the structure to be Type7; if $k<0$, then we predict the structure to be Type5. With the three steps above, we complete the prediction of structure in the approximate solver. The algorithm for predicting the structure is summarised in Algorithm\ref{algorithm: prediction}.
\begin{algorithm}[t]
	\caption{Structure prediction algorithm}
	\label{algorithm: prediction}
	\hspace*{0.02in} {\bf Input:}
	left initial value $U_L$, right initial value $U_R$,
	coefficients $k_1,k_2,k_3,\gamma$ \\
	\hspace*{0.02in} {\bf Output:} 
	The structure of Riemann solution
	\begin{algorithmic}[1]
		\If{($k>0$ and $M_L\geq M_2^*$) or ($k=0$ and $M_L\geq 1$) or ($k<0$ and $1\leq M_L<M_3^{**}$)}
		\If{$p< p_+$}
		\State \Return Type2
		\ElsIf{$(\gamma+1)p\leq (2\gamma M_+^2-\gamma+1)p_+$}
		\State \Return Type2
		\EndIf
		\ElsIf{$T(p_1)T(p_2)<0$}
		\State \Return Type1
		\Else
		\If{$k>0$}
		\State \Return Type3
		\ElsIf{$k=0$}
		\State \Return Type7
		\Else
		\State \Return Type5
		\EndIf
		\EndIf
	\end{algorithmic}
\end{algorithm}

\textbf{Step2: solutions of Riemann solver:}

For different predicted structures, we use different methods to obtain approximate solutions $\widetilde{U}_-$ and $\widetilde{U}_+$, as shown in Table\ref{table: approximate solutions}.
\begin{table}[htbp]\label{table: approximate solutions}
	\centering
	\caption{The approximate solutions for each sturcture}
	\begin{tabular}{ccc}
		\toprule
		Structure & left-hand approximate solution $\widetilde{U}_-$& right-hand approximate solution $\widetilde{U}_+$\\
		\midrule
		\specialrule{0em}{3pt}{3pt}
		Type1& $W(U_L)\cap \{U|p=T^{-1}(0)\}$&\multirow{4}{*}{$D(\widetilde{U}_-,k_1,k_2,k_3)$}\\	
		\specialrule{0em}{1pt}{1pt}	
		Type2& $U_L$&\\		
		\specialrule{0em}{1pt}{1pt}	
		Type3& $W(U_L)\cap \{U|M=M_1^*\}$&\\	
		\specialrule{0em}{1pt}{1pt}		
		Type5 and Type7& $R(U_L)\cap \{U|M=1\}$&\\
		\specialrule{0em}{3pt}{3pt}
		\bottomrule
	\end{tabular}
\end{table}

The branch used for the solution $D(\widetilde{U}_-,k_1,k_2,k_3)$ of $\widetilde{U}_+$ is specified by the structure. The only thing that needs to be given is the solution of $T^{-1}(0)$.
The root of the function $T(p)$ is obtained by a dichotomy in our program. 
The initial values of the iteration are shown below.
\begin{equation}\label{solver: initial value of iteration}
	\begin{cases}
		\text{if } p_2<p_L<p_1 \text{ and } T(p_L)*T(p_1)\leq 0, \text{ then the initial value of iteration are } p_1 \text{ and } p_L,\\
		\text{if } p_2<p_L<p_1 \text{ and } T(p_L)*T(p_2)\leq 0, \text{ then the initial value of iteration are } p_2 \text{ and } p_L,\\
		\text{else}, \text{ then the initial value of iteration are } p_1 \text{ and } p_2.
	\end{cases}
\end{equation}

The following conclusion shows that the proposed solver is able to give an exact solution when the initial value of the Riemann problem is a single stationary wave.
\begin{Lemma}\label{lemma: solver}
	If the initial value of Riemann problem satisfies $M_L\in \Gamma_-$ and 
	$U_R\in \mathcal{D}(U_L,k_1,k_2,k_3)$, then the solution of the approximate solver is $\widetilde{U}_-=U_L, \widetilde{U}_+=U_R$.
\end{Lemma}
\begin{Proof}
	We divide the proof into two steps, corresponding to the cases where the initial value is a supersonic stationary wave and a subsonic stationary wave, respectively.
	
	We first show that the conclusion holds when the initial value is a supersonic stationary wave($M_L,M_R\geq 1$). Obviously, $U_L$ satisfies the condition (\ref{solver: type2 prediction1}). Since $U_+=U_R$, the condition (\ref{solver: type2 prediction2}) can also be satisfied, therefore the predicted structure is Type2.
	According to Table\ref{table: approximate solutions}, the solution of solver is $\widetilde{U}_-=U_L, \widetilde{U}_+=U_R$.
	
	Secondly, we show that the conclusion holds when the initial value is a subsonic stationary wave. If $k\geq 0$, then we have $M_L<1,M_R\leq 1$, and the condition (\ref{solver: type2 prediction1}) cannot be satisfied. If $k<0$, then we have $M_L\leq 1,M_R<1$, and it is easy to check that the condition (\ref{solver: type2 prediction2}) cannot be satisfied. Thus the predicted structure cannot be Type2. By
	\begin{equation*}
		0<f_M^{1-W}(U_L,p_L)=M_L\leq \begin{cases}
			M_*^1,\ k>0\\
			1,\ k\leq 0
		\end{cases}
	\end{equation*}
	and the monotonicity of $f_M^{1-W}$, we have $p_2\leq p_L<p_1$. 
	\begin{equation*}
		\begin{aligned}
			T(p_L)=&f_u^{3-\mathcal{W}}(U_R,p_Lg_p^\mathcal{D}(f_M^{1-\mathcal{W}}(U_L,p_L)))-f_u^{1-\mathcal{W}}(U_L,p_L)g_u^\mathcal{D}(f_M^{1-\mathcal{W}}(U_L,p_L))\\
			=&f_u^{3-\mathcal{W}}(U_R,p_Lg_p^\mathcal{D}(M_L))-u_Lg_u^\mathcal{D}(M_L)
			=f_u^{3-\mathcal{W}}(U_R,p_R)-u_R=0\\
		\end{aligned}
	\end{equation*}
	Because $p_2\leq p_L<p_1$ and $p_L$ is a root of $T(p)$, the solution to the dichotomy with the initial value (\ref{solver: initial value of iteration}) is $U_L$, i.e. the solution to the approximate solver is $\widetilde{U}_-=U_L, \widetilde{U}_+=U_R$.
\end{Proof}

\subsubsection{solver-based numerical flux}
We apply the approximate Riemann solver above to the construction of numerical fluxes:
\begin{equation}\label{new scheme: flux}
	\widehat{F}_{{j_0}+\frac{1}{2},-}^n=F(\widetilde{U}_-),
	\widehat{F}_{{j_0}+\frac{1}{2},+}^n=F(\widetilde{U}_+),
\end{equation}
where $\widetilde{U}_-$ and $\widetilde{U}_+$ are given by the approximate Riemann solver with the initial values $U(x_{{j_0}+\frac{1}{2}}^-,t)$ and $U(x_{{j_0}+\frac{1}{2}}^+,t)$.

It follows from Lemma\ref{lemma: solver} that the numerical fluxes (\ref{new scheme: flux}) satisfy the consistency condition (\ref{consistency condition 2}).

\begin{Theorem}\label{theorem: solver-based scheme}
	Under the grid (\ref{grid}), if the numerical fluxes at the cell interfaces other than the origin satisfy the consistency condition (\ref{consistency condition 1}), and the numerical fluxes at the origin are (\ref{new scheme: flux}), then the conservative scheme (\ref{unsplitting conservative scheme}) is well-balanced. 
\end{Theorem}

In the scheme based on (\ref{new scheme: flux}), the branches of stationary wave curve used to calculate $\widehat{F}_{{j_0}+\frac{1}{2},-}^n$ and $\widehat{F}_{{j_0}+\frac{1}{2},-}^n$ are chosen based on the structure of the Riemann solution, so that they are always the same. 
In addition, due to the consideration of the coupling of convection and source terms, our scheme not only maintains the equilibrium states, but also gives satisfactory results for unsteady solutions that contain multiple elementary waves.

\section{Numerical tests}\label{Numerical tests}

A series of Riemann problems with different structures are designed to test the effects of these three numerical schemes. All of these tests use $\gamma=1.4$. The coefficients of source terms, initial values, and structures for each test are shown in Table\ref{table:test}.

\begin{table}[htbp]
	\centering
	\caption{The coefficients, initial values and structures of Test 1$\sim$8.}
	\label{table:test}
	\begin{tabular}{ccccc}
		\toprule
		\specialrule{0em}{3pt}{3pt}
		\multirow{2}*{Test}&coefficients&$U_L$&$U_R$&\multirow{2}*{Structure}\\
		\specialrule{0em}{1pt}{1pt}
		&$(k_1,k_2,k_3)$& $(\rho_L,u_L,p_L)$&$(\rho_R,u_R,p_R)$&\\
		\specialrule{0em}{3pt}{3pt}	
		\midrule
		\specialrule{0em}{3pt}{3pt}
		1&$(0.4,0.2,0.4)$& $(0.6,0.5,0.6)$&$(0.641338,0.654881,0.62495)$&{single source stationary wave}\\		
		\specialrule{0em}{1pt}{2pt}
		2&$(0.2,0.0,0.2)$& $(1.0,1.0,1.0)$&$(0.933943,0.411564,1.27555)$&Type1\\
		\specialrule{0em}{1pt}{2pt}
		3&$(0.1,0.1,0.2)$& $(1.0,2.0,1.0)$&$(1.888,2.53245,2.17219)$&Type2\\
		\specialrule{0em}{1pt}{2pt}
		4&$(0.1,-0.2,0.2)$& $(1.0,1.0,1.0)$&$(0.378535,2.07562,0.46455)$&Type3\\
		\specialrule{0em}{1pt}{2pt}
		5&$(0.1,0.1,0.2)$& $(1.0,1.74007,1.0)$&$(1.65217,1.71963,1.77224)$&Type4\\
		\specialrule{0em}{1pt}{2pt}
		6&$(0.1,0.2,-0.2)$& $(1.0,0.8,1.0)$&$(1.27959,1.38758,0.671459)$&Type5\\
		\specialrule{0em}{1pt}{2pt}
		7&$(0.1,0.1,-0.1)$& $(1.0,0.8,1.0)$&$(0.468365,1.92334,0.450956)$&Type6\\
		\specialrule{0em}{1pt}{2pt}
		8&$(0.3,0.3,0.3)$& $(0.6,0.8,0.6)$&$(0.459223,1.45488,0.507773)$&Type7\\
		\specialrule{0em}{3pt}{3pt}
		\bottomrule
	\end{tabular}
\end{table}

We use a 3rd-order discontinuous Galerkin(DG) method\cite{1989TVB} as a spatial discretization method.
A characteristicwise TVD (total variation diminishing) limiter is empolyed to avoid numerical oscillations (see \cite{1989TVB}).
For the discretization of the convection equation (\ref{splitting: PDE}) in the splitting method we use a local Lax-Friedrichs numerical flux:
\begin{equation}
	\widehat{F}^{LF}(U_1,U_2)=\frac{1}{2}\left[ F(U_1)+F(U_2)-\alpha (U_2-U_1) \right] , 
\end{equation}
where $\alpha =\max_{1\leq p\leq 3}(|\lambda^{(p)}(U_1)|,|\lambda^{(p)}(U_2)|)$.

The local Lax-Friedrichs numerical flux above is also used
for the numerical fluxes (\ref{unsplitting: flux other than origin}) in the unsplitting scheme and
the approximation to the classical Riemann solution (\ref{W-B: flux}) in the K-T scheme.

A third-order Runge-Kutta method \cite{Gottlieb2001Strong} is used for time discretization. The wave induced by the source term is a stationary discontinuity, and the influence of source term propagates through convection, hence a standard Courant-Friedrichs-Levy (CFL) condition is sufficient:
\begin{equation*}
	\sup\limits_{i,x}\left| \lambda_i(U^h(x,t^0)) \right| 
	\frac{\Delta t}{h}\leq \frac{1}{2}.
\end{equation*}

\begin{figure}[htbp]	\label{test1-1}
	\centering
	\subfigure[density]{
		\includegraphics[width=0.45\linewidth]{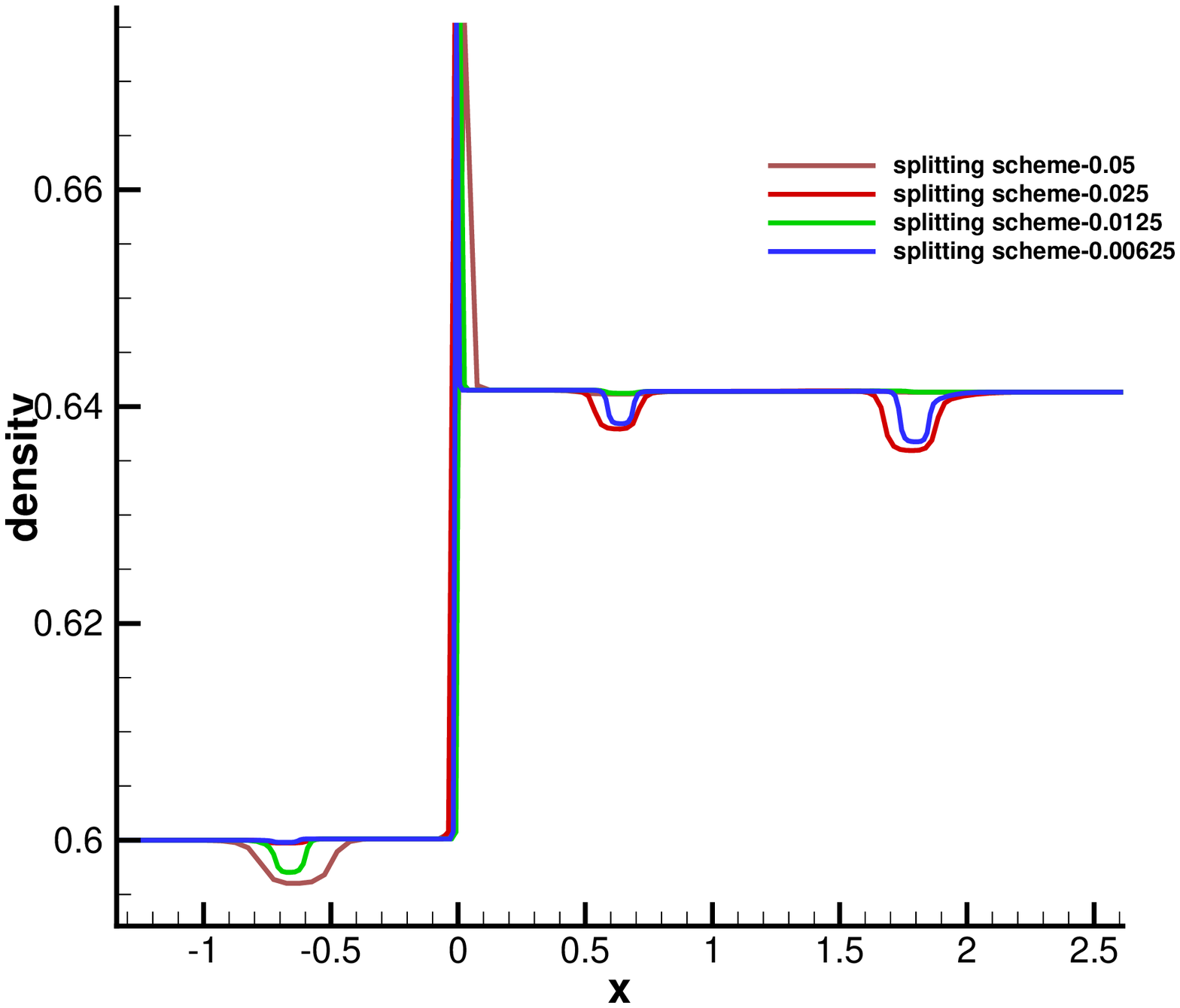}
	}
	\quad
	\subfigure[pressure]{
		\includegraphics[width=0.45\linewidth]{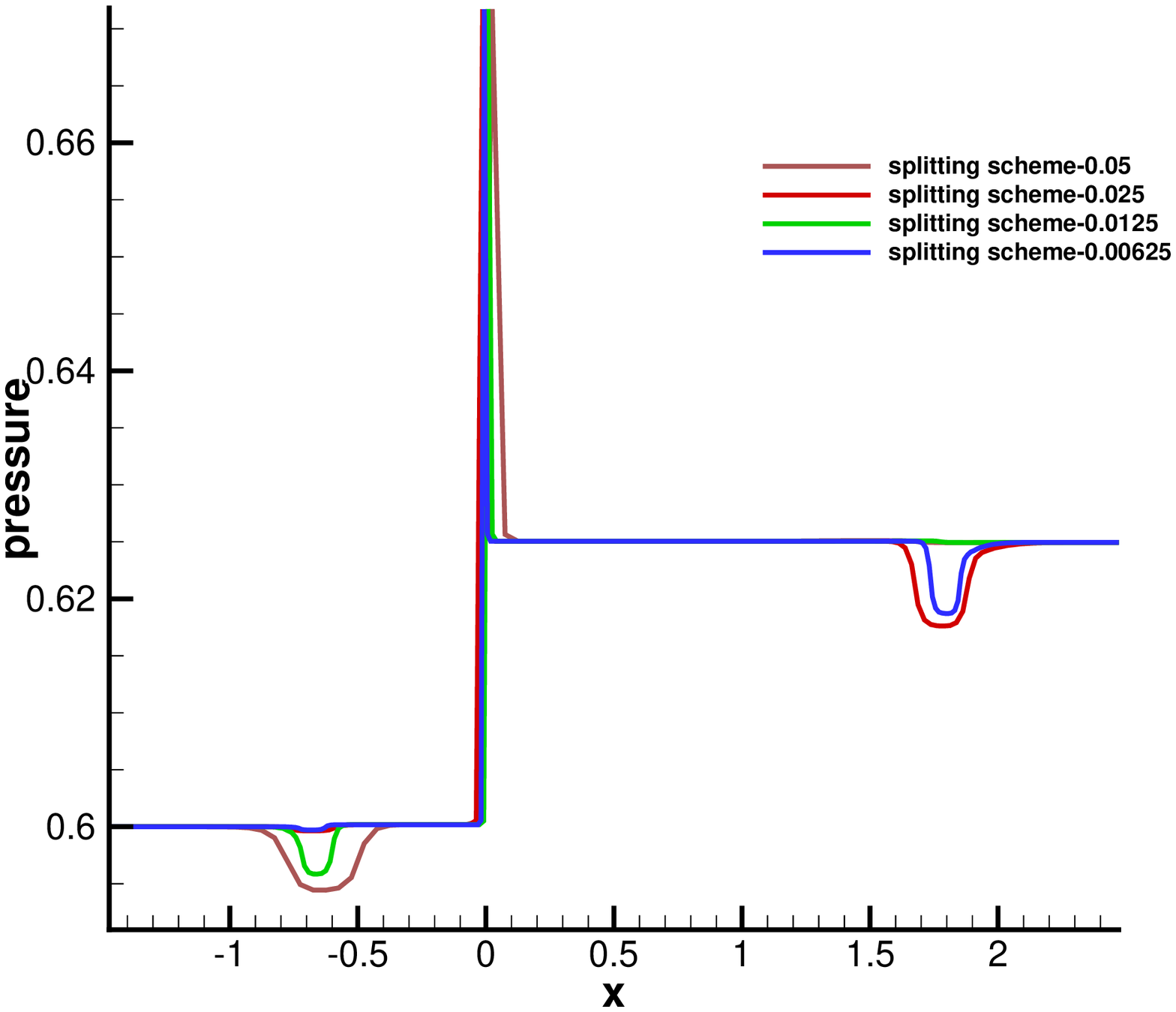}
	}
	\caption{The numerical solution of splitting scheme of Test1 at time $t=1.0s$ with different $h$: 0.05, 0.025, 0.0125, 0.00625. Left: density. Right: pressure.}
\end{figure}

\begin{figure}[htbp]
	\centering
	\begin{minipage}[c]{0.49\textwidth}
		\centering		
		\includegraphics[width=0.9\linewidth]{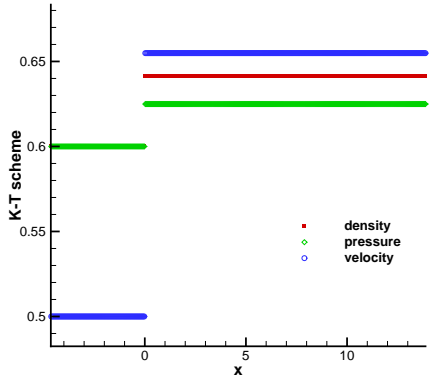}
		\caption{The numerical solution of K-T scheme of Test1 with $h=0.05$ at time $t=1.0s$.}
		\label{test1-2}
	\end{minipage}
	\begin{minipage}[c]{0.49\textwidth}		
		\centering		
		\includegraphics[width=0.9\linewidth]{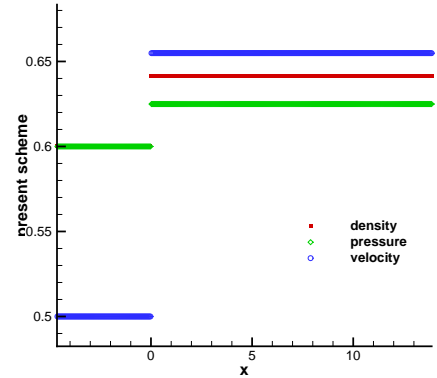}
		\caption{The numerical solution of present scheme of Test1 with $h=0.05$ at time $t=1.0s$.}
		\label{test1-3}
	\end{minipage}
\end{figure}

\begin{figure}	
	\centering	
	\begin{minipage}[c]{0.45\textwidth}		
		\centering		
		\subfigure[density]{
			\includegraphics[width=0.9\linewidth]{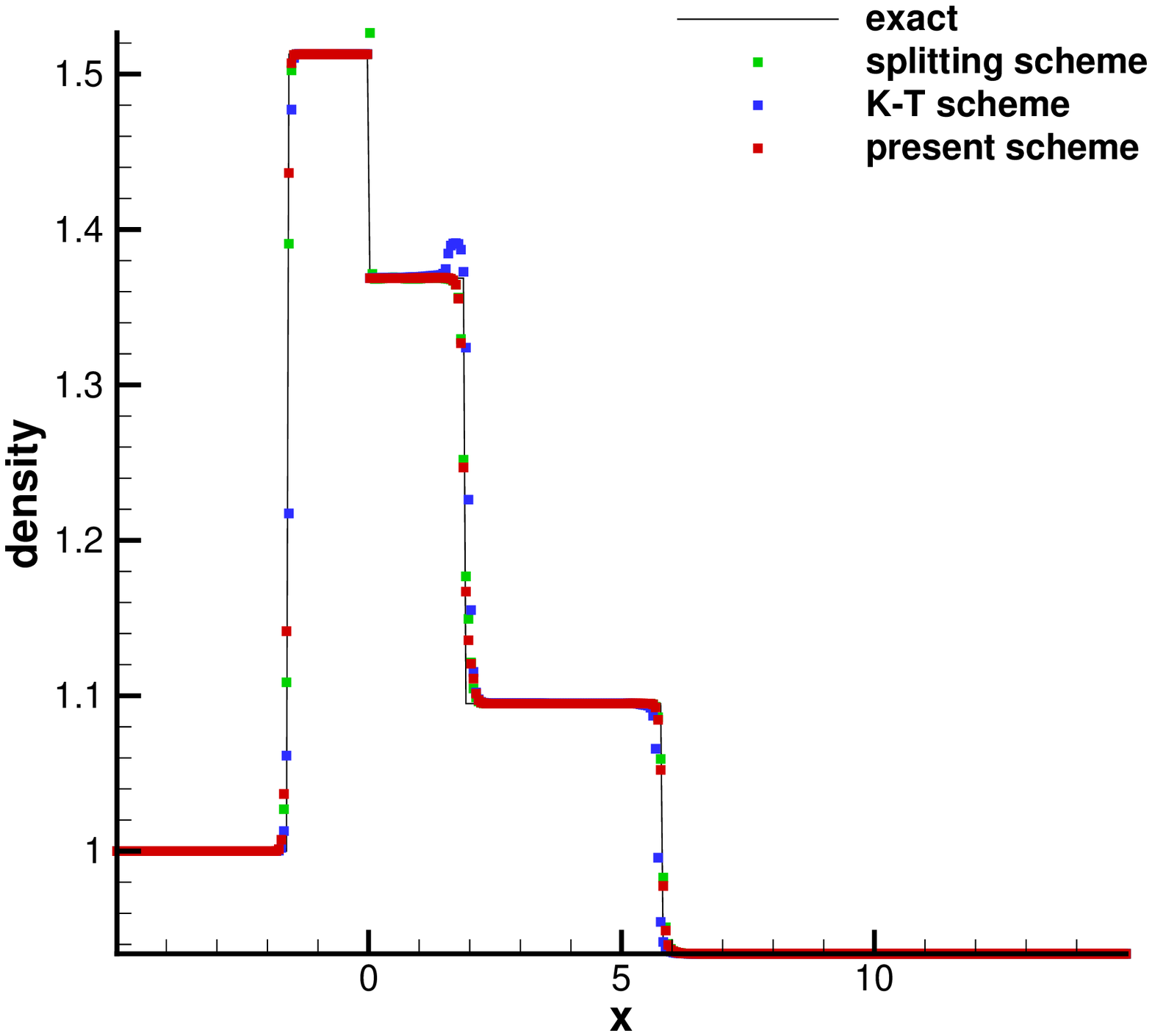}
		}
		\\
		\subfigure[velocity]{
			\includegraphics[width=0.9\linewidth]{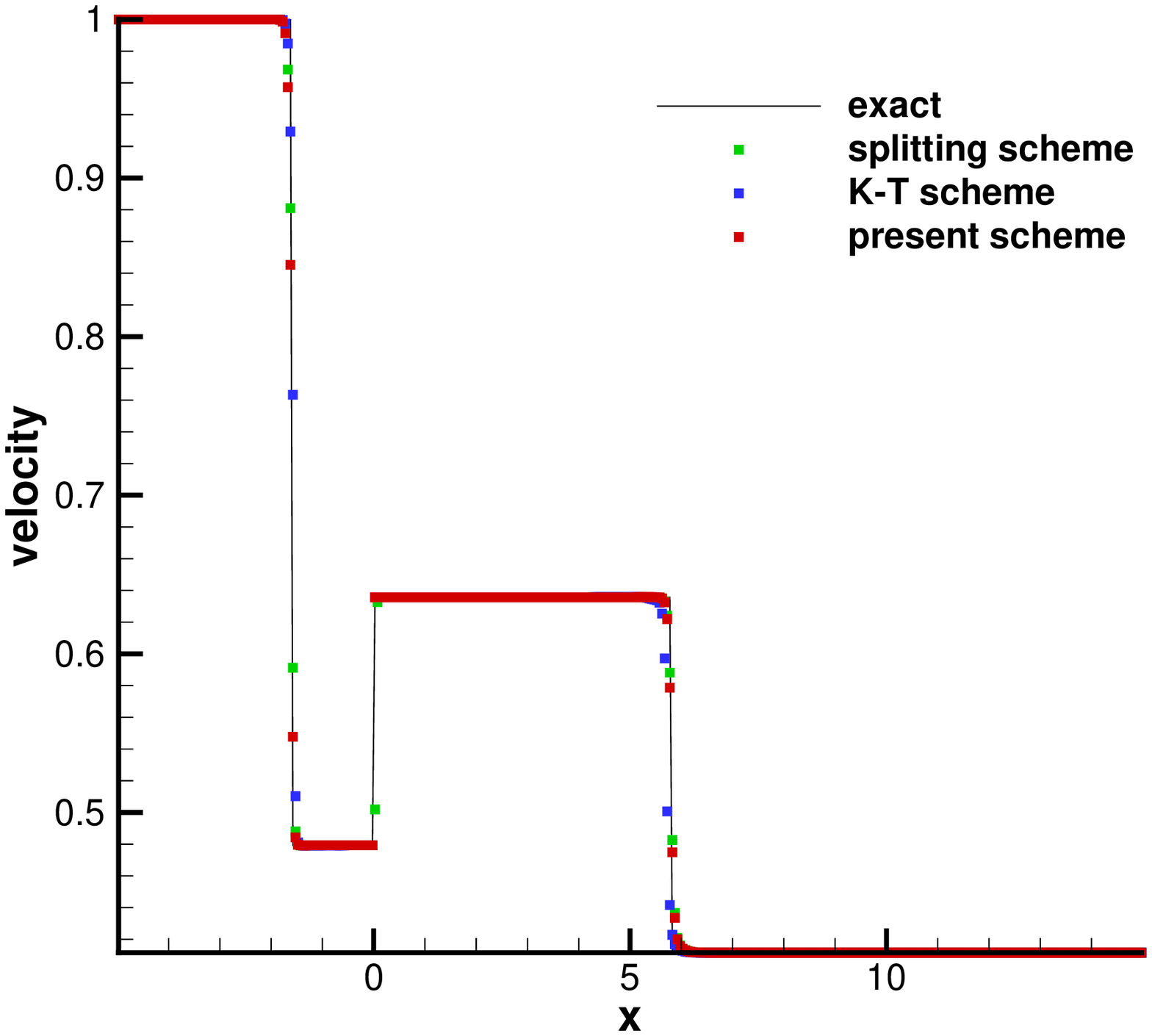}
		}
		\\
		\subfigure[pressure]{
			\includegraphics[width=0.9\linewidth]{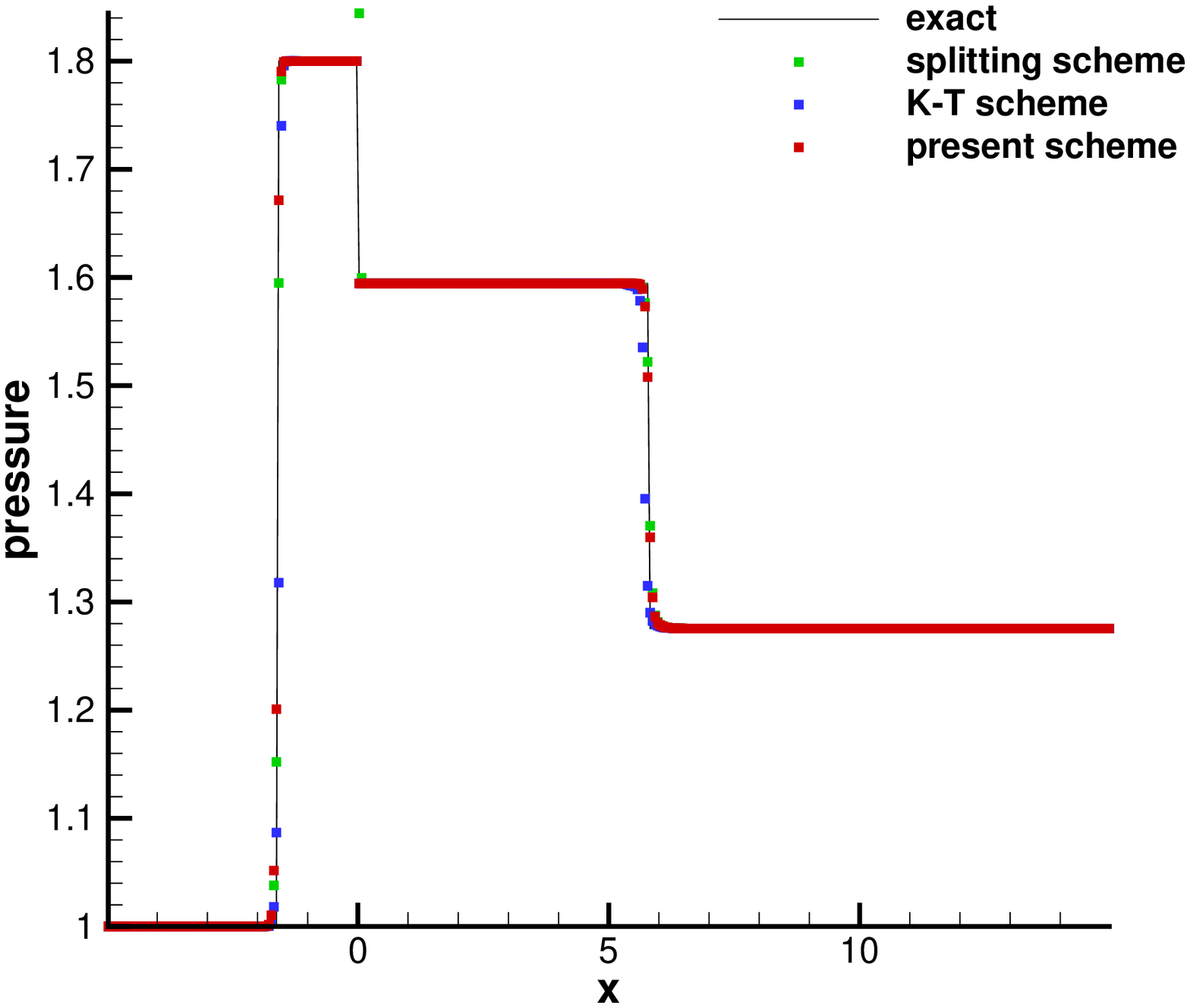}
		}
	\end{minipage}%
	\begin{minipage}[c]{0.45\textwidth}				
		\centering		
		\subfigure[density]{
			\includegraphics[width=0.9\linewidth]{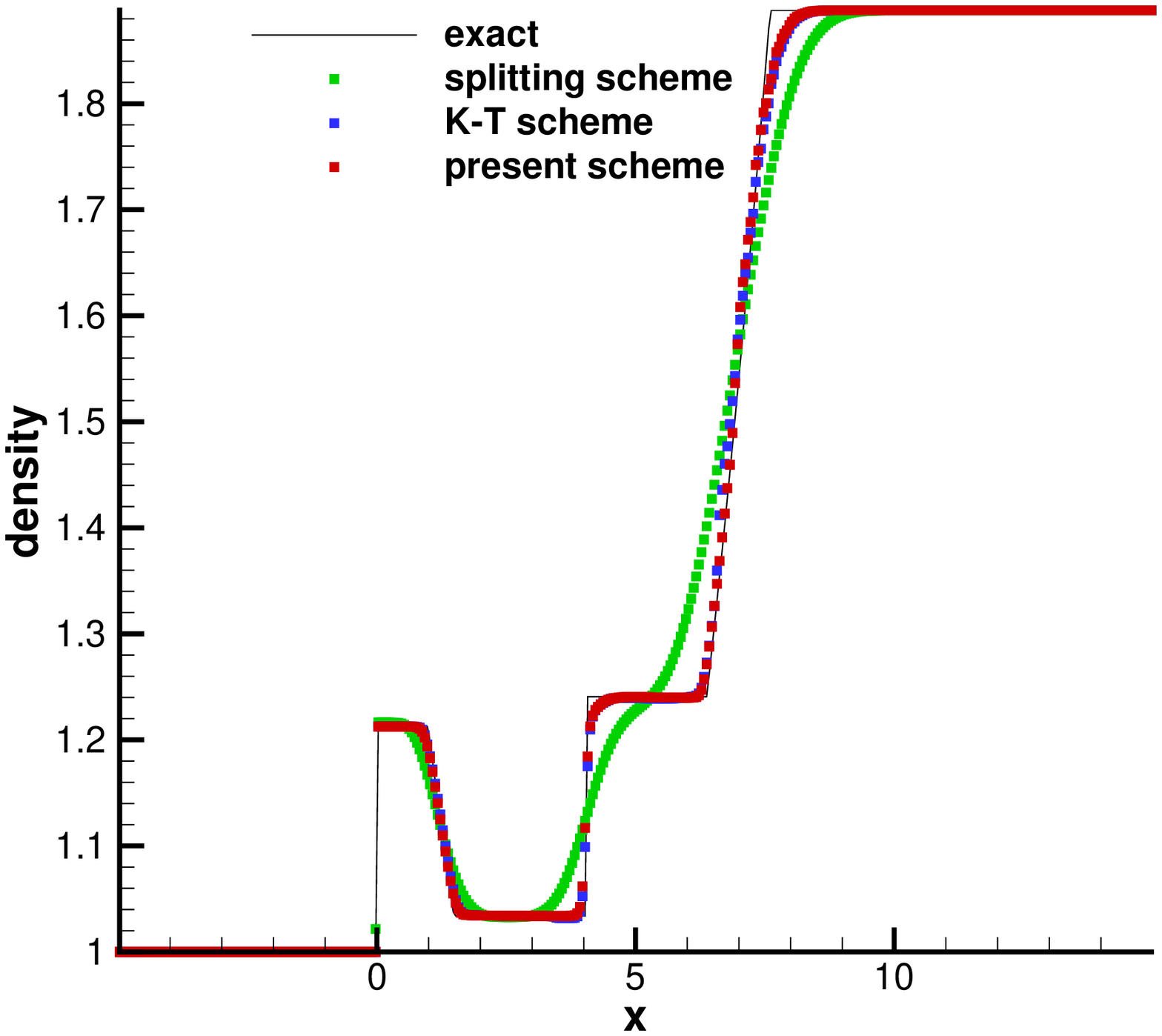}
		}
		\\
		\subfigure[velocity]{
			\includegraphics[width=0.9\linewidth]{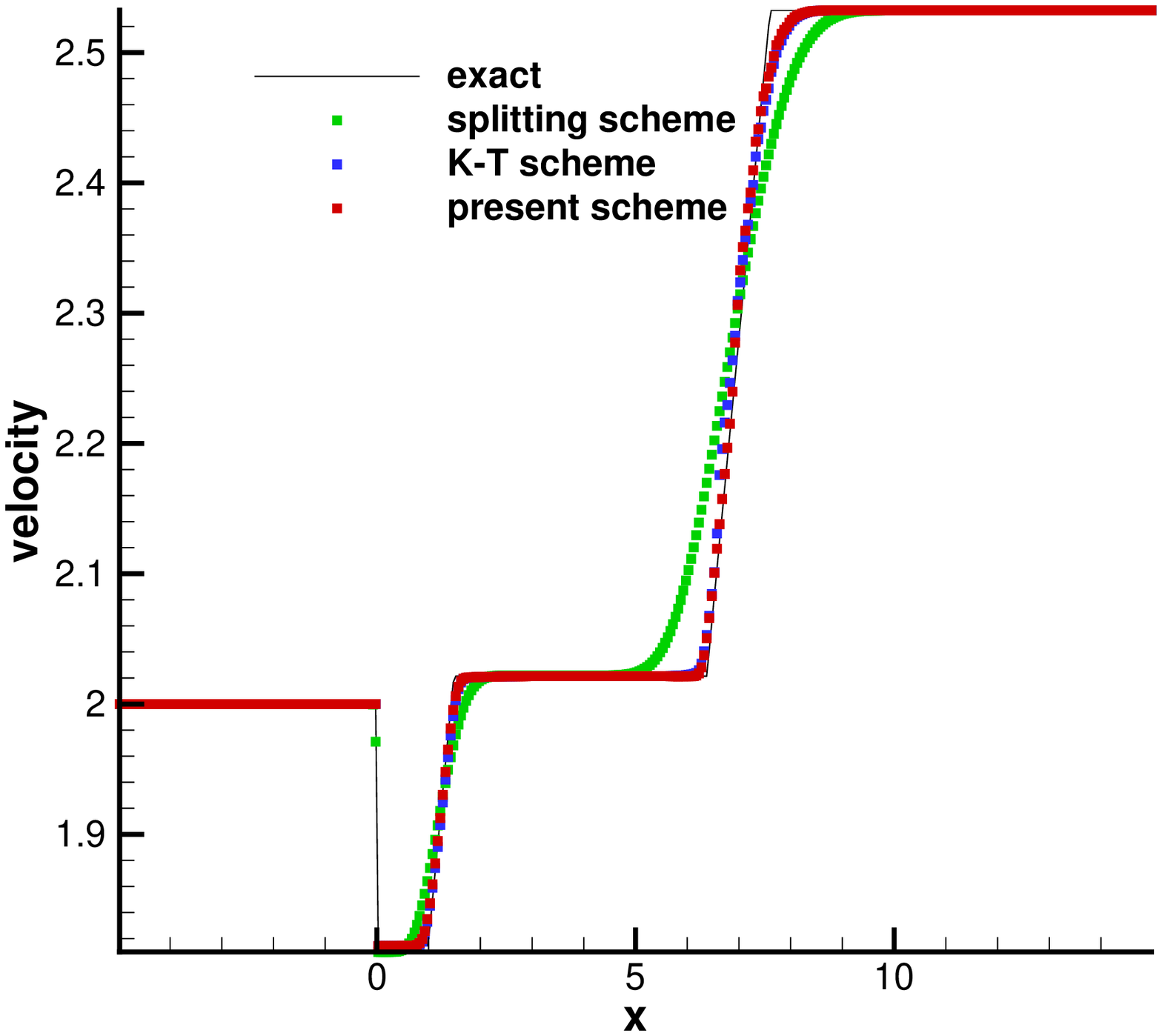}
		}
		\\
		\subfigure[pressure]{
			\includegraphics[width=0.9\linewidth]{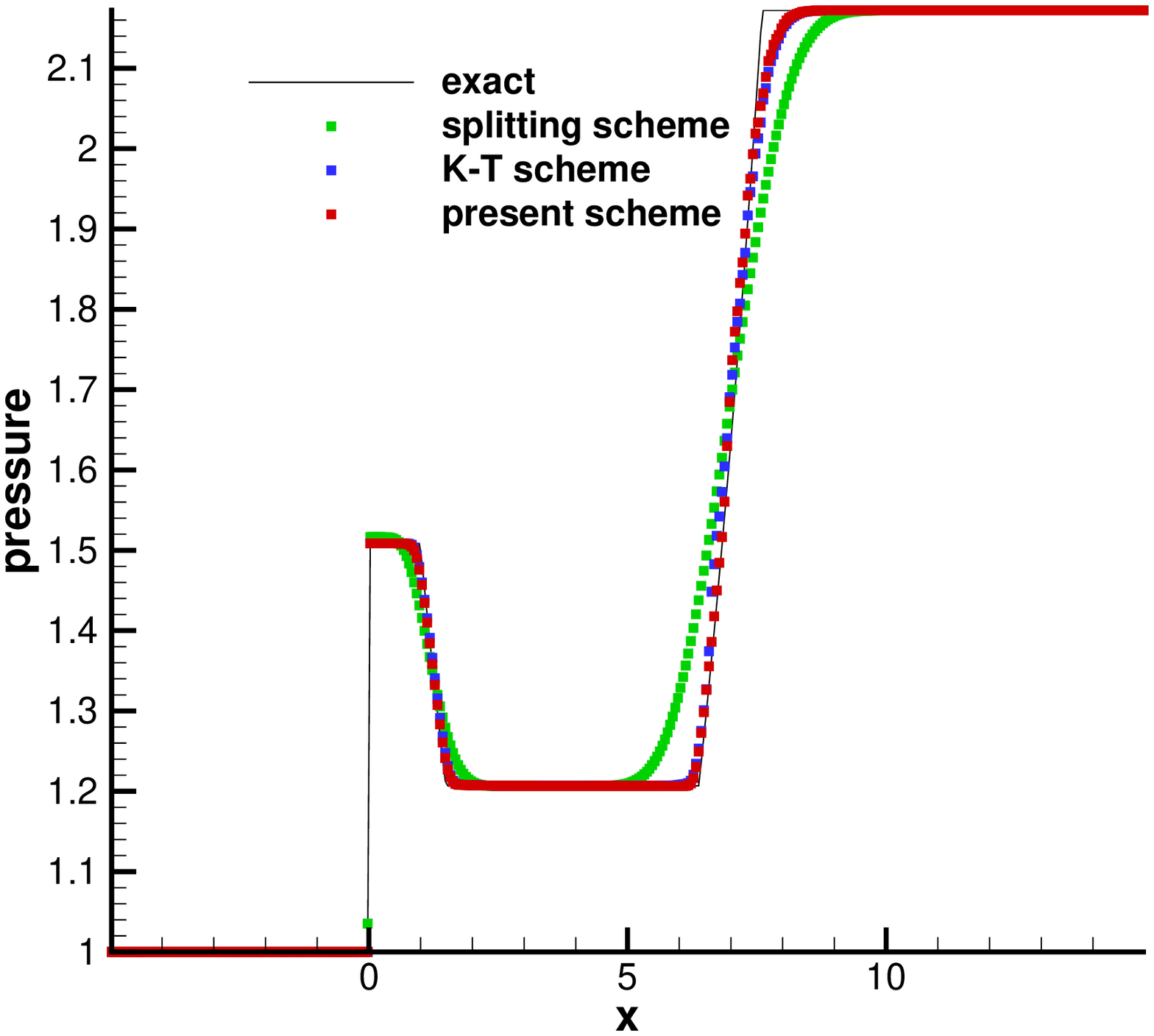}
		}	
	\end{minipage}	
	\caption{Splitting scheme(green symbols), K-T scheme(blue symbols) and present scheme(red symbols) versus exact solutions(solid line) with $h=0.05$. Left column: Test2 and $t=3.0s$, right column: Test3 and $t=2.0s$. Top: density, middle: velocity, and bottom: pressure.}
	\label{test2-3}
\end{figure}

\begin{figure}	
	\centering	
	\begin{minipage}[c]{0.5\textwidth}				
		\centering		
		\subfigure[density]{
			\includegraphics[width=0.9\linewidth]{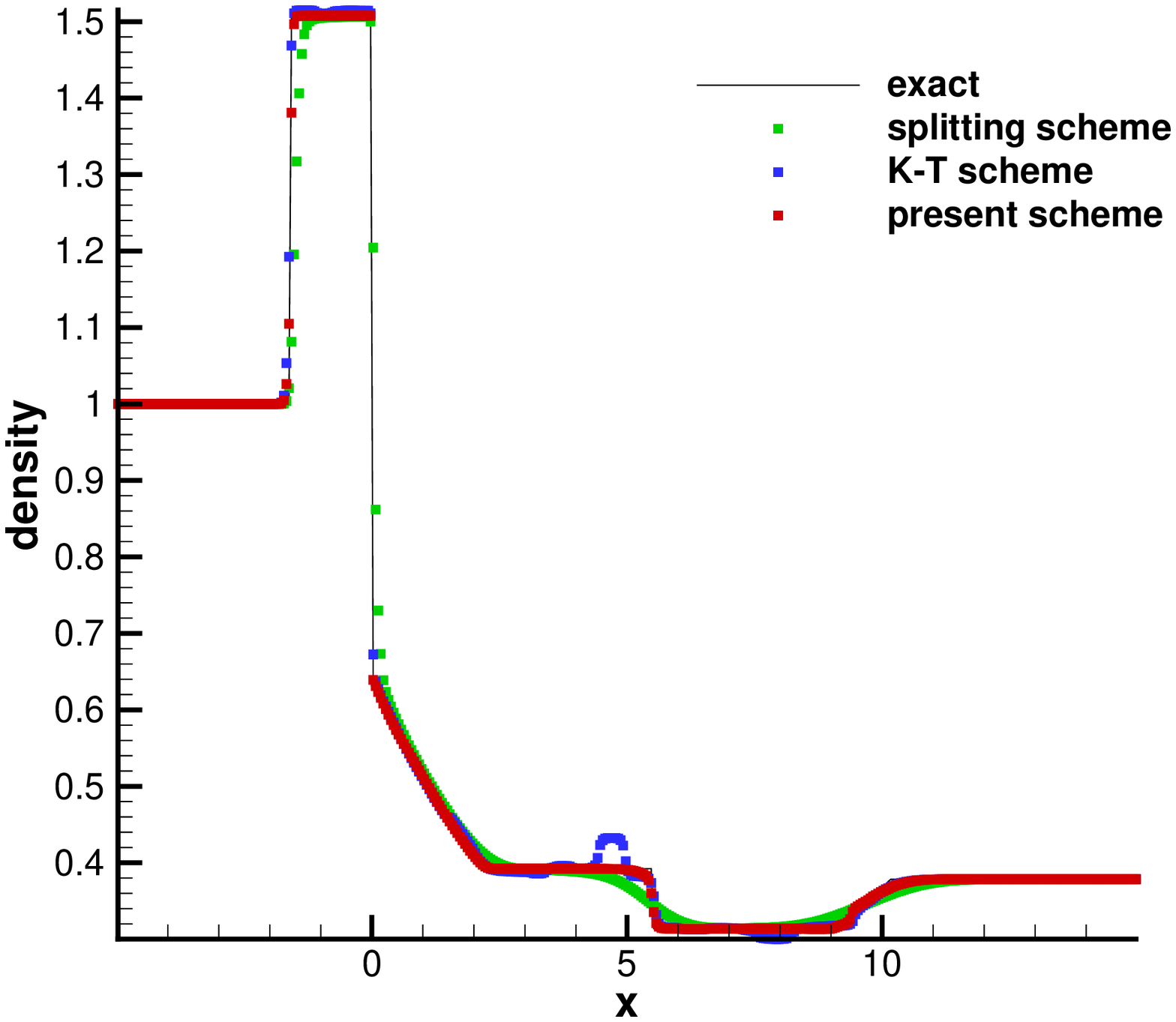}
		}
		\\
		\subfigure[velocity]{
			\includegraphics[width=0.9\linewidth]{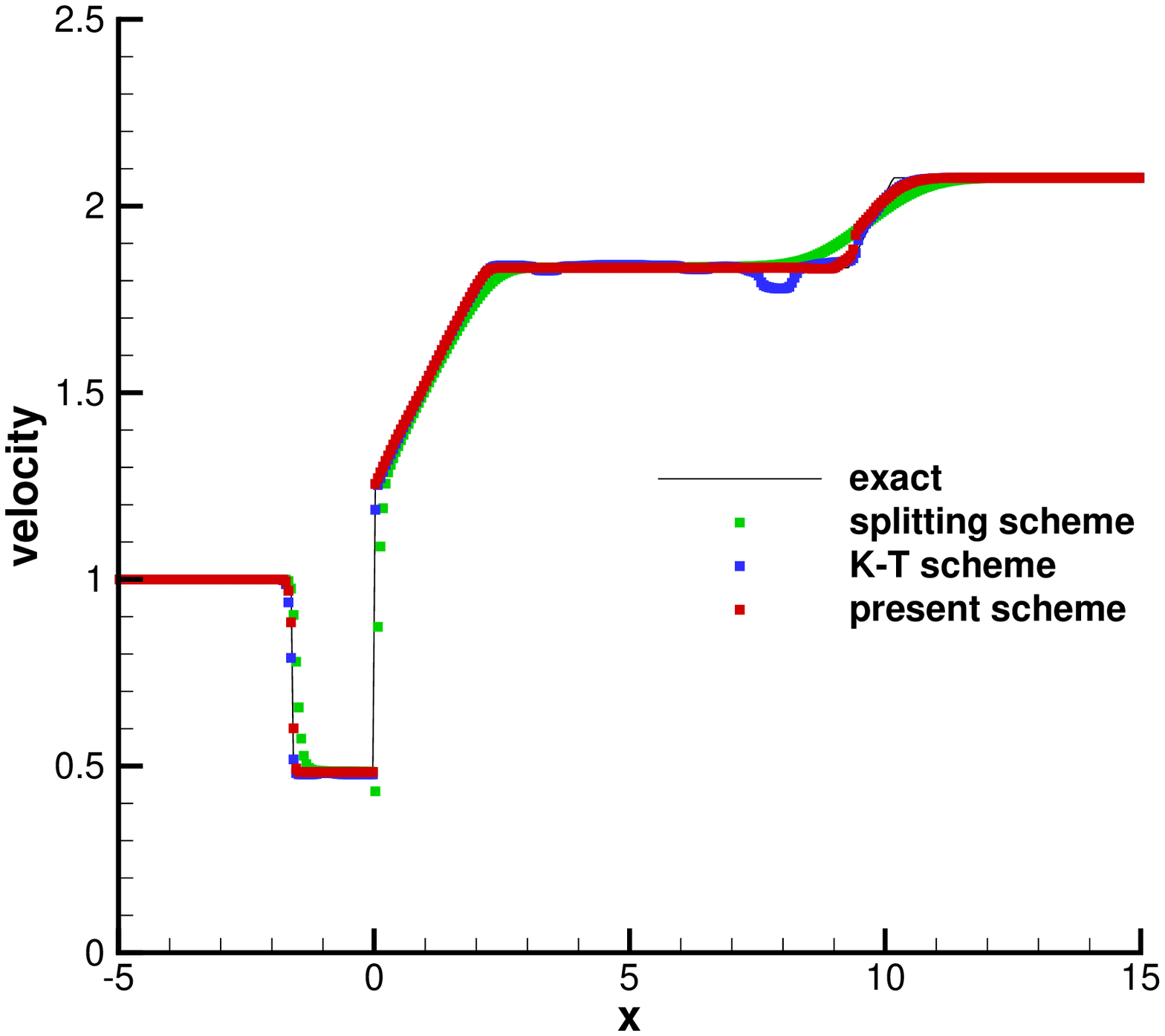}
		}
		\\
		\subfigure[pressure]{
			\includegraphics[width=0.9\linewidth]{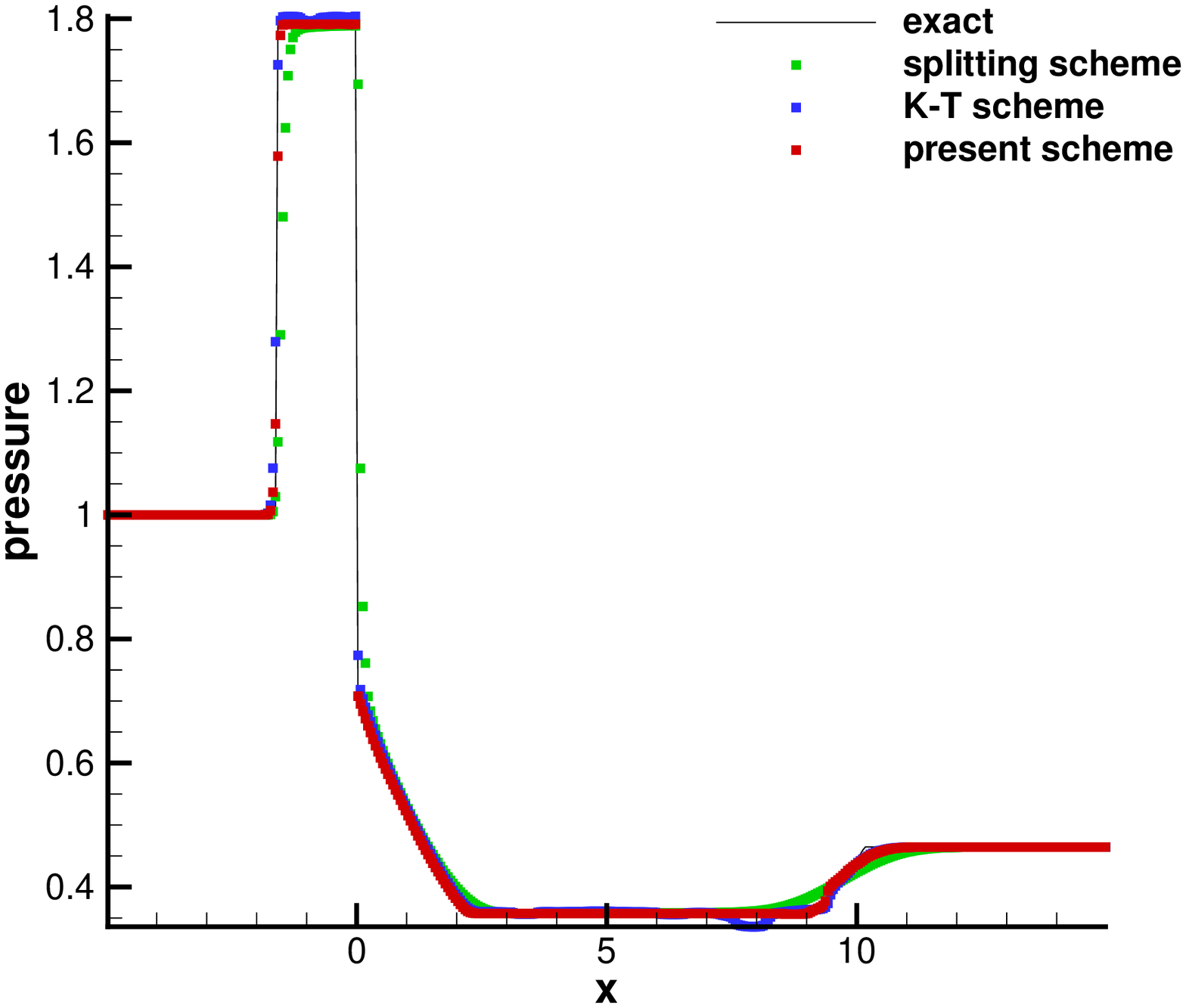}
		}	
	\end{minipage}%
	\begin{minipage}[c]{0.5\textwidth}				
		\centering		
		\subfigure[density]{
			\includegraphics[width=0.9\linewidth]{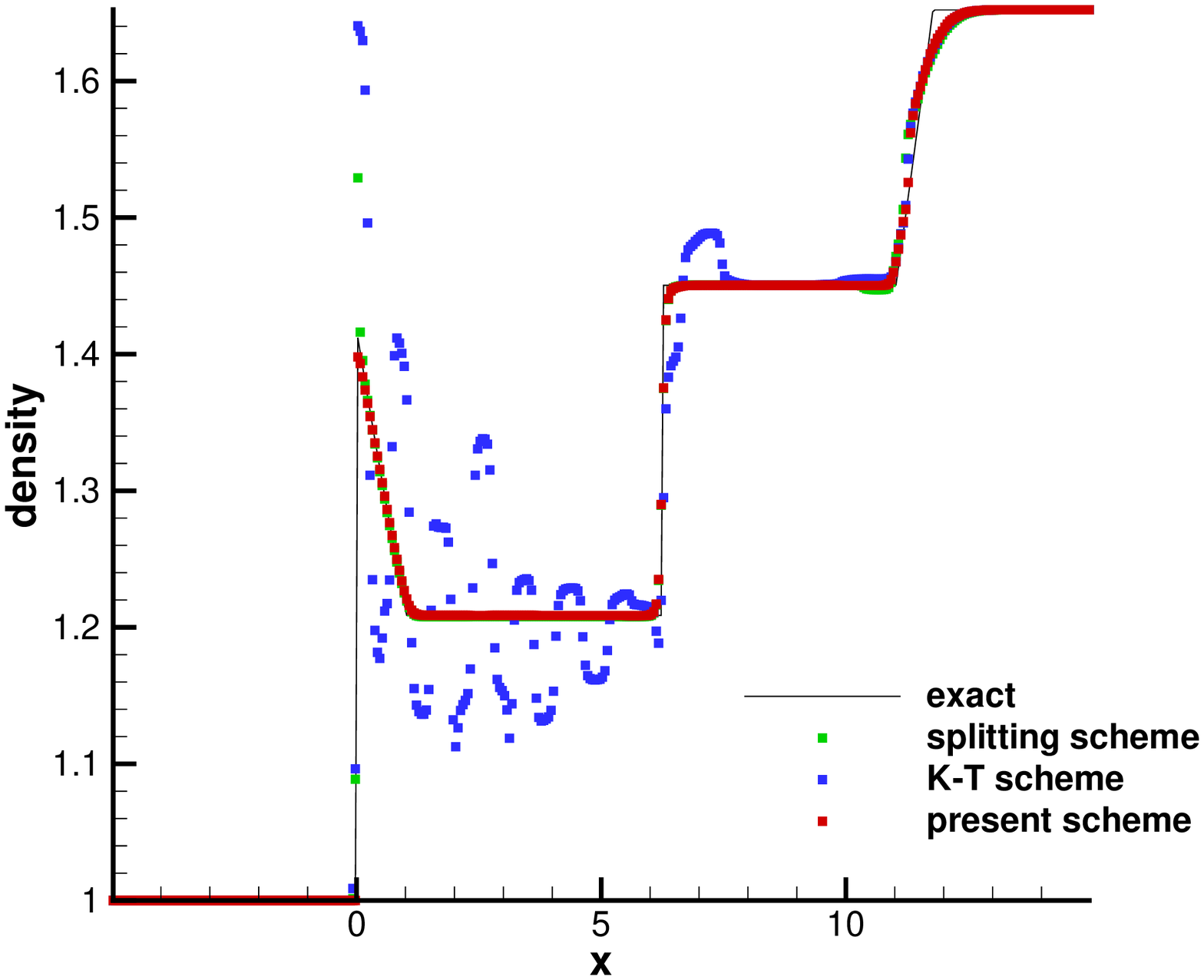}
		}
		\\
		\subfigure[velocity]{
			\includegraphics[width=0.9\linewidth]{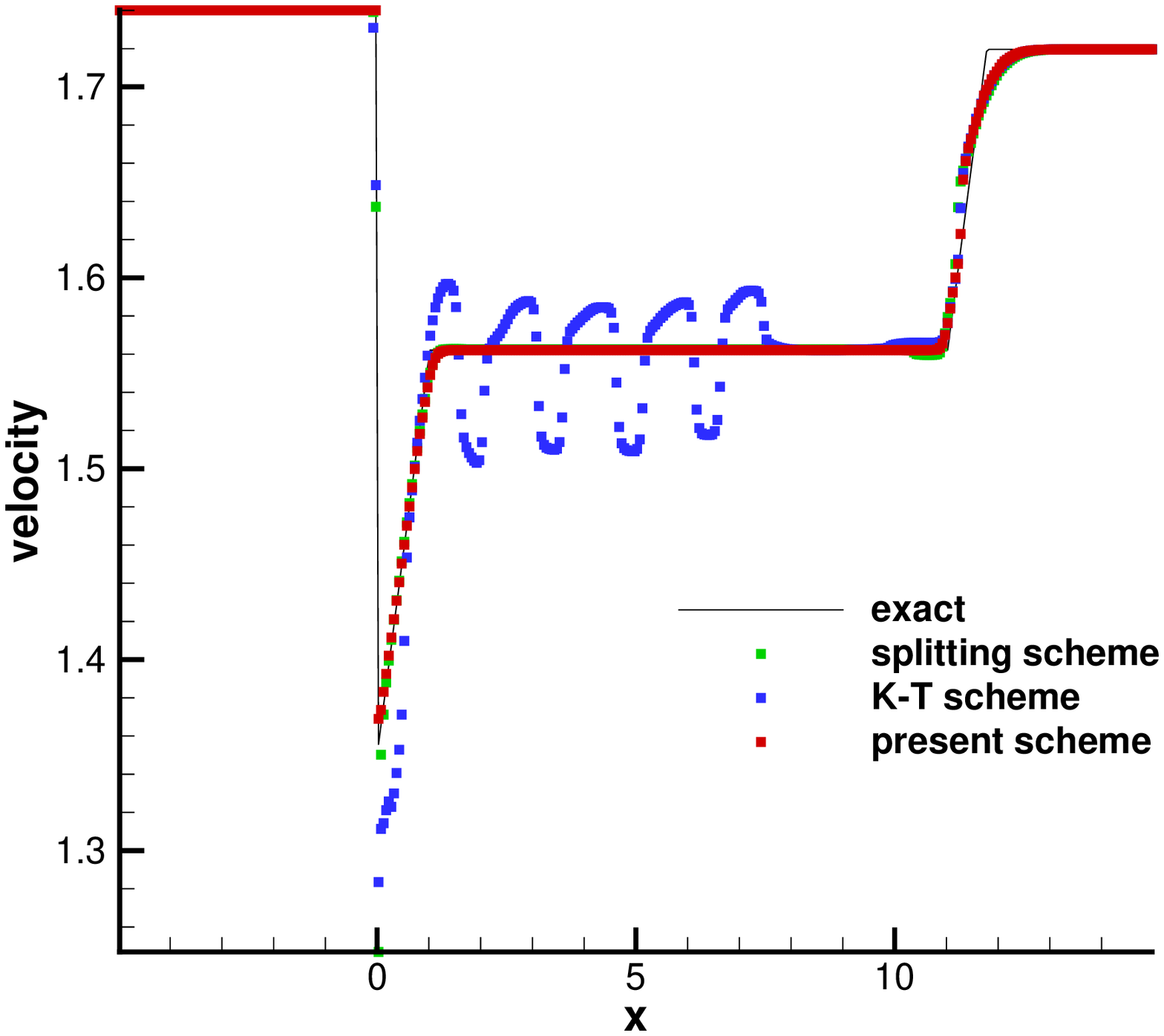}
		}
		\\
		\subfigure[pressure]{
			\includegraphics[width=0.9\linewidth]{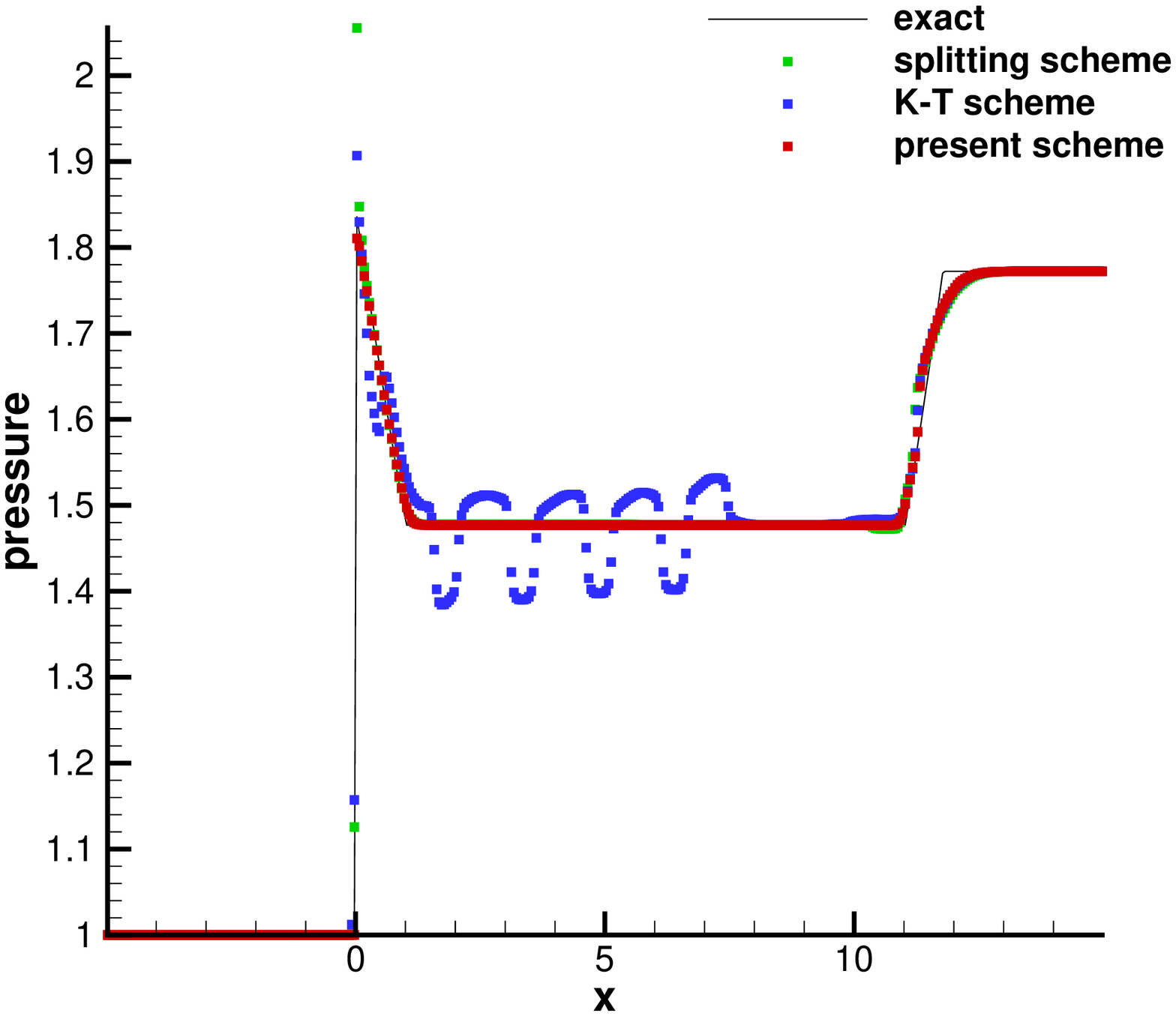}
		}
	\end{minipage}	
	\caption{Splitting scheme(green symbols), K-T scheme(blue symbols) and present scheme(red symbols) versus exact solutions(solid line) with $h=0.05$. Left column: Test4 and $t=3.0s$, right column: Test5 and $t=4.0s$. Top: density, middle: velocity, and bottom: pressure.}
	\label{test4-5}
\end{figure}

\begin{figure}	
	\centering	
	\begin{minipage}[c]{0.5\textwidth}				
		\centering		
		\subfigure[density]{
			\includegraphics[width=0.9\linewidth]{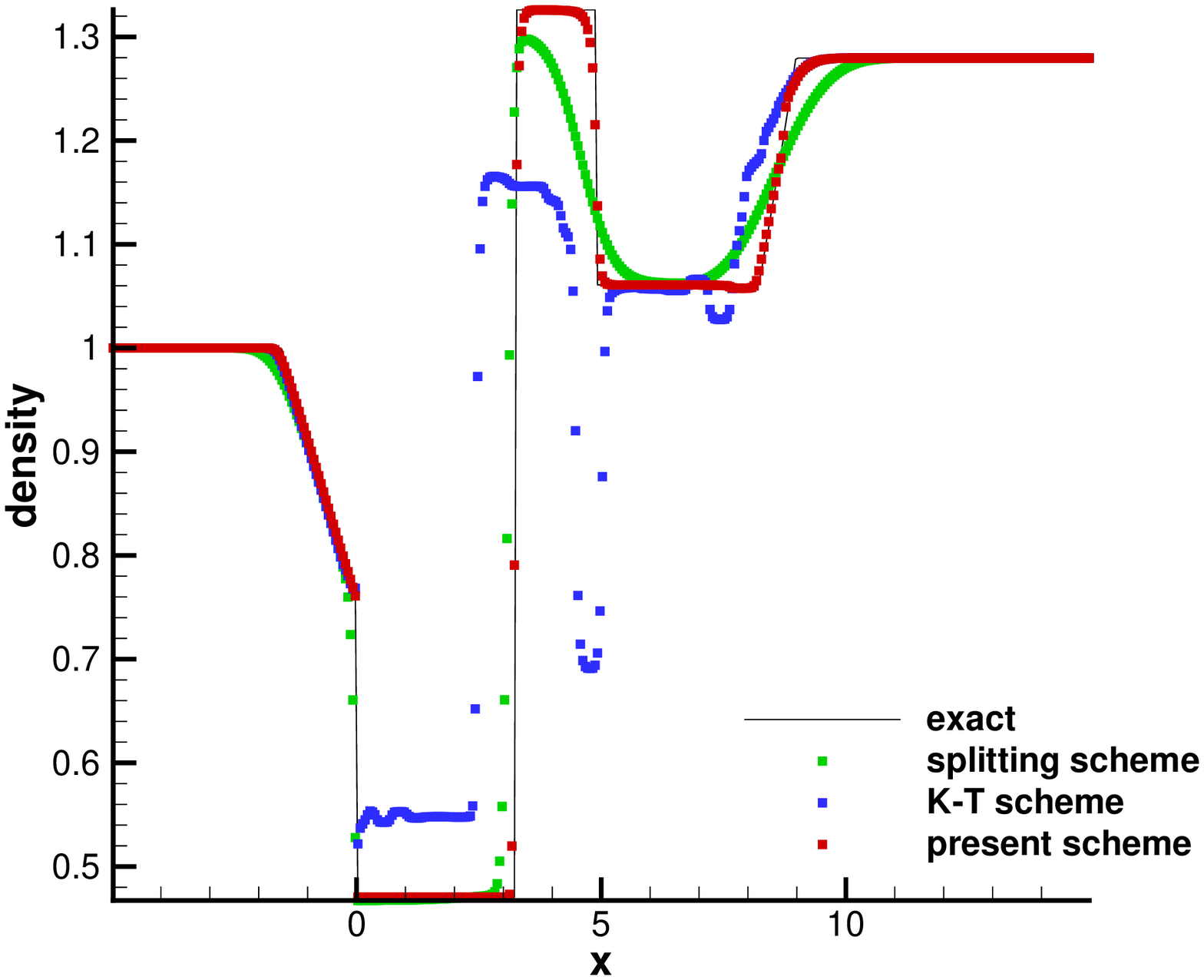}
		}
		\\
		\subfigure[velocity]{
			\includegraphics[width=0.9\linewidth]{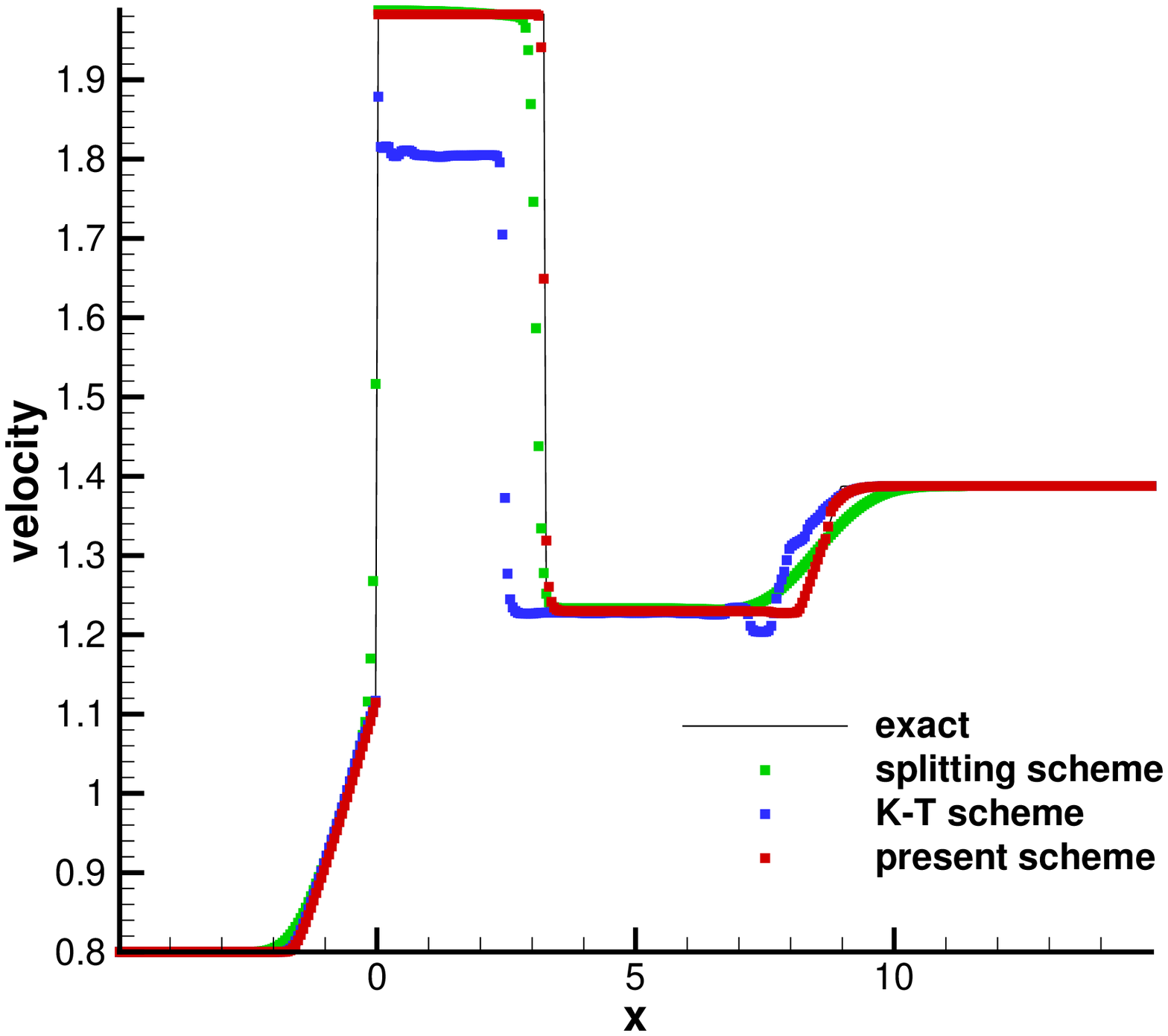}
		}
		\\
		\subfigure[pressure]{
			\includegraphics[width=0.9\linewidth]{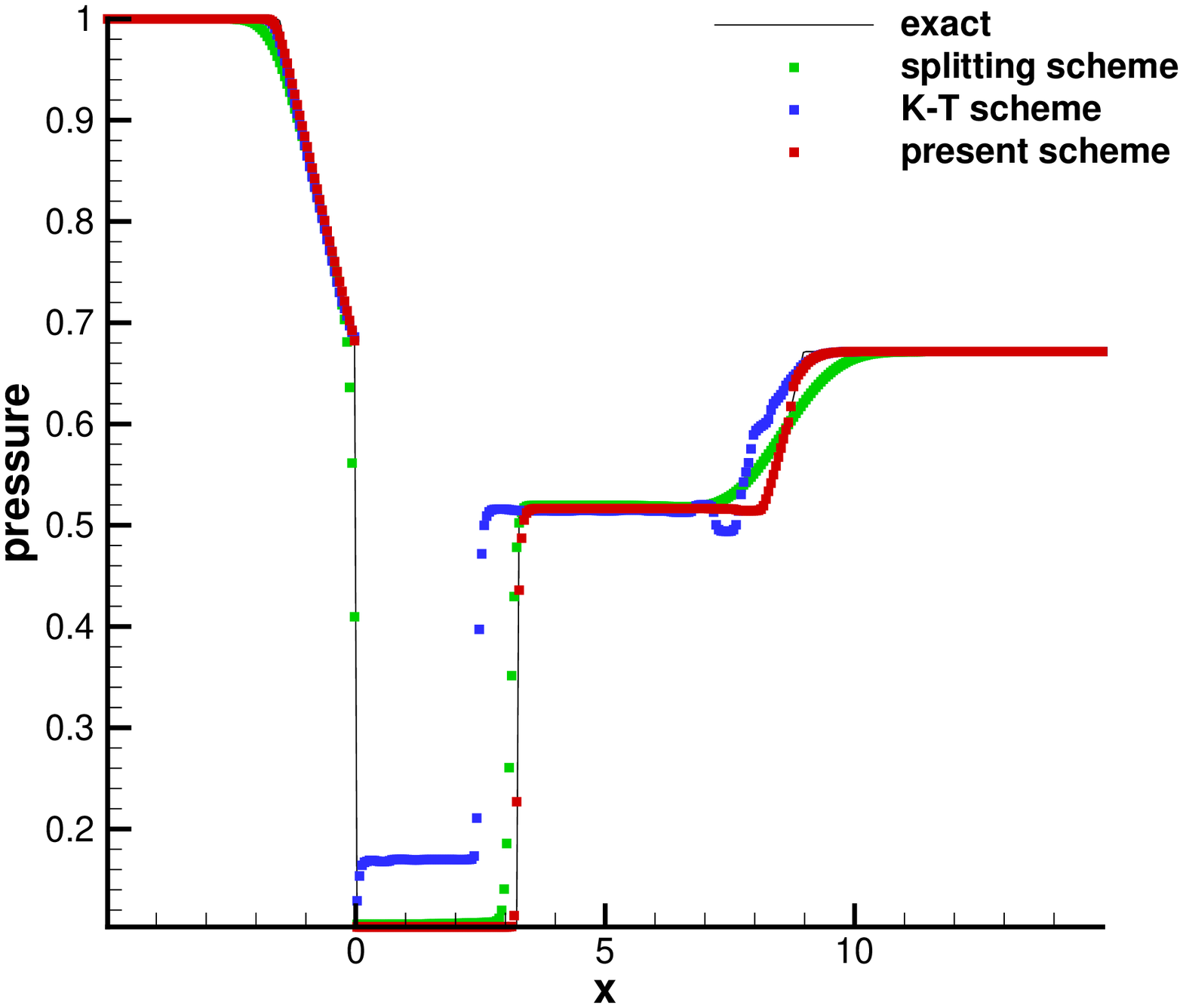}
		}	
	\end{minipage}%
	\begin{minipage}[c]{0.5\textwidth}				
		\centering		
		\subfigure[density]{
			\includegraphics[width=0.9\linewidth]{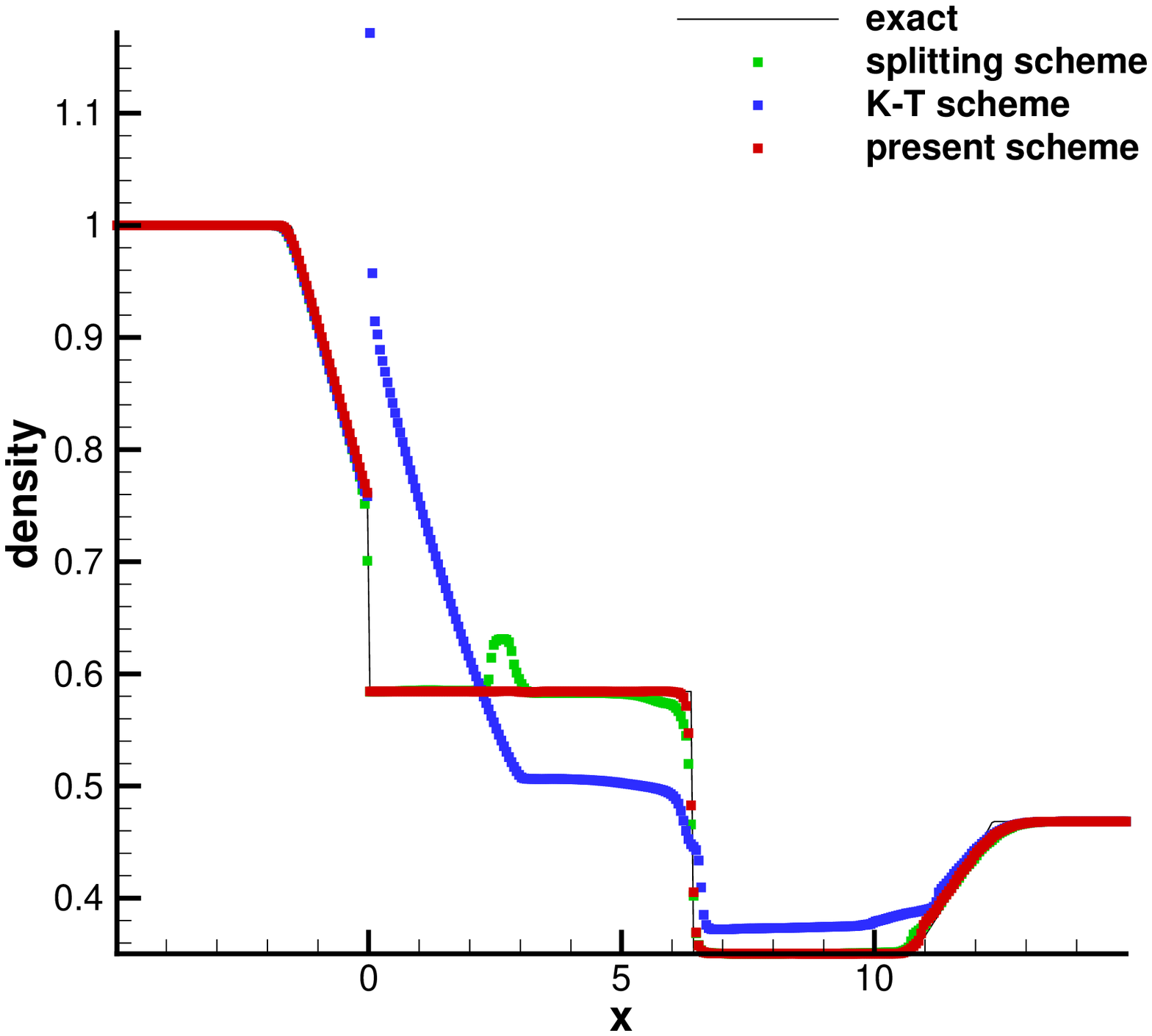}
		}
		\\
		\subfigure[velocity]{
			\includegraphics[width=0.9\linewidth]{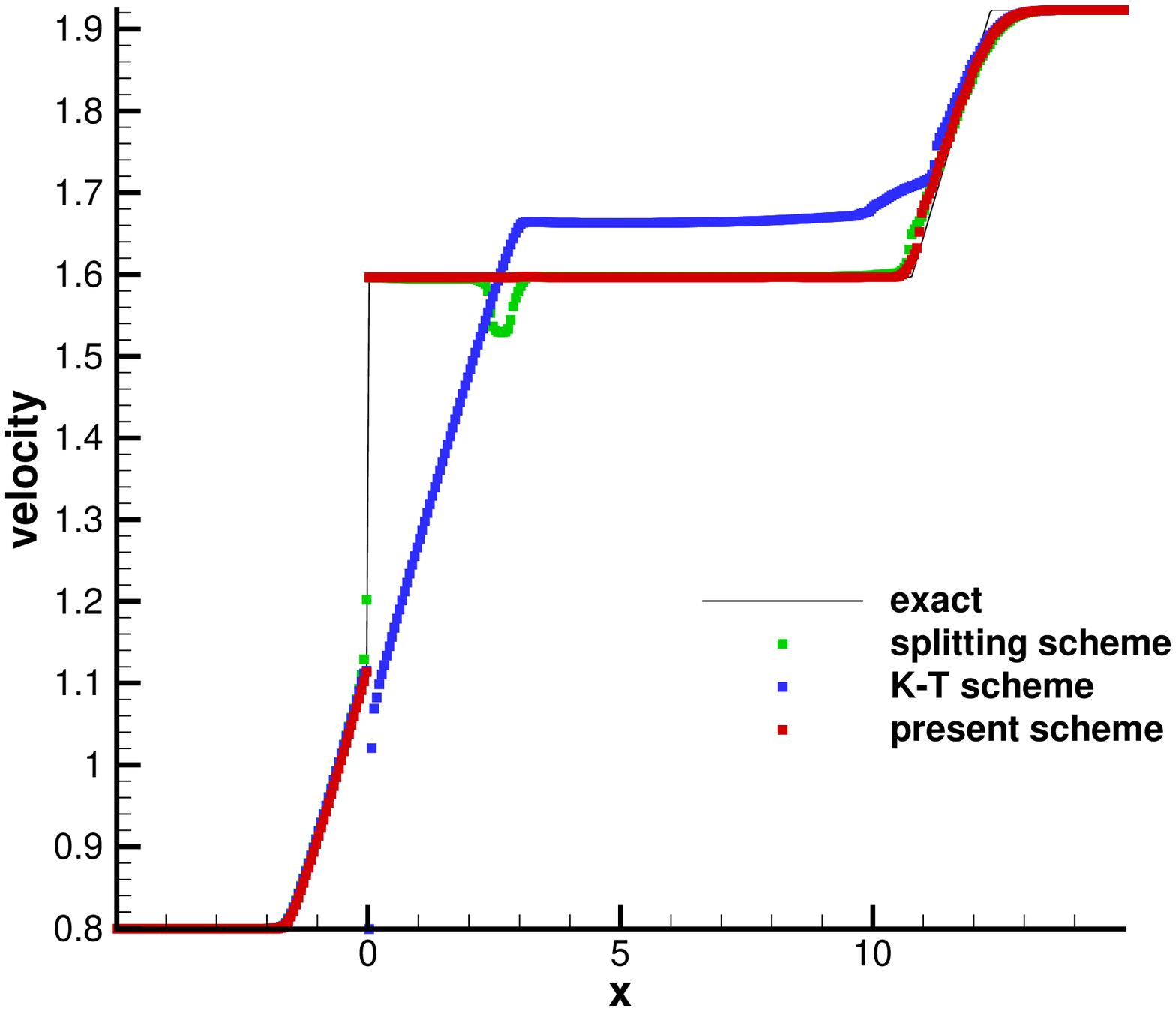}
		}
		\\
		\subfigure[pressure]{
			\includegraphics[width=0.9\linewidth]{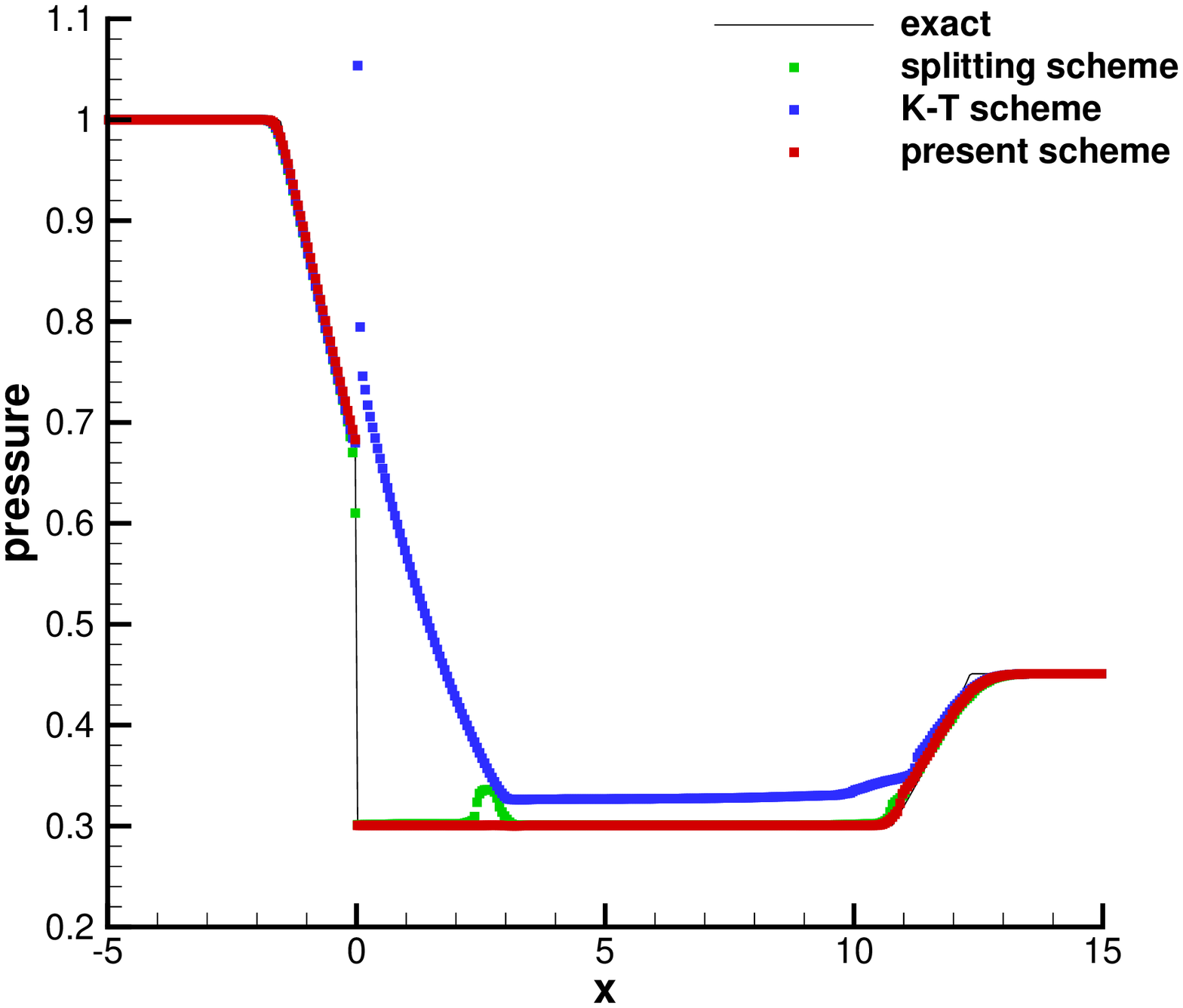}
		}
	\end{minipage}	
	\caption{Splitting scheme(green symbols), K-T scheme(blue symbols) and present scheme(red symbols) versus exact solutions(solid line) with $h=0.05$. Left column: Test6 and $t=4.0s$, right column: Test7 and $t=4.0s$. Top: density, middle: velocity, and bottom: pressure.}
	\label{test6-7}
\end{figure}

\begin{figure}[htbp]		
	\centering
	\subfigure[density]{
		\includegraphics[width=0.4\linewidth]{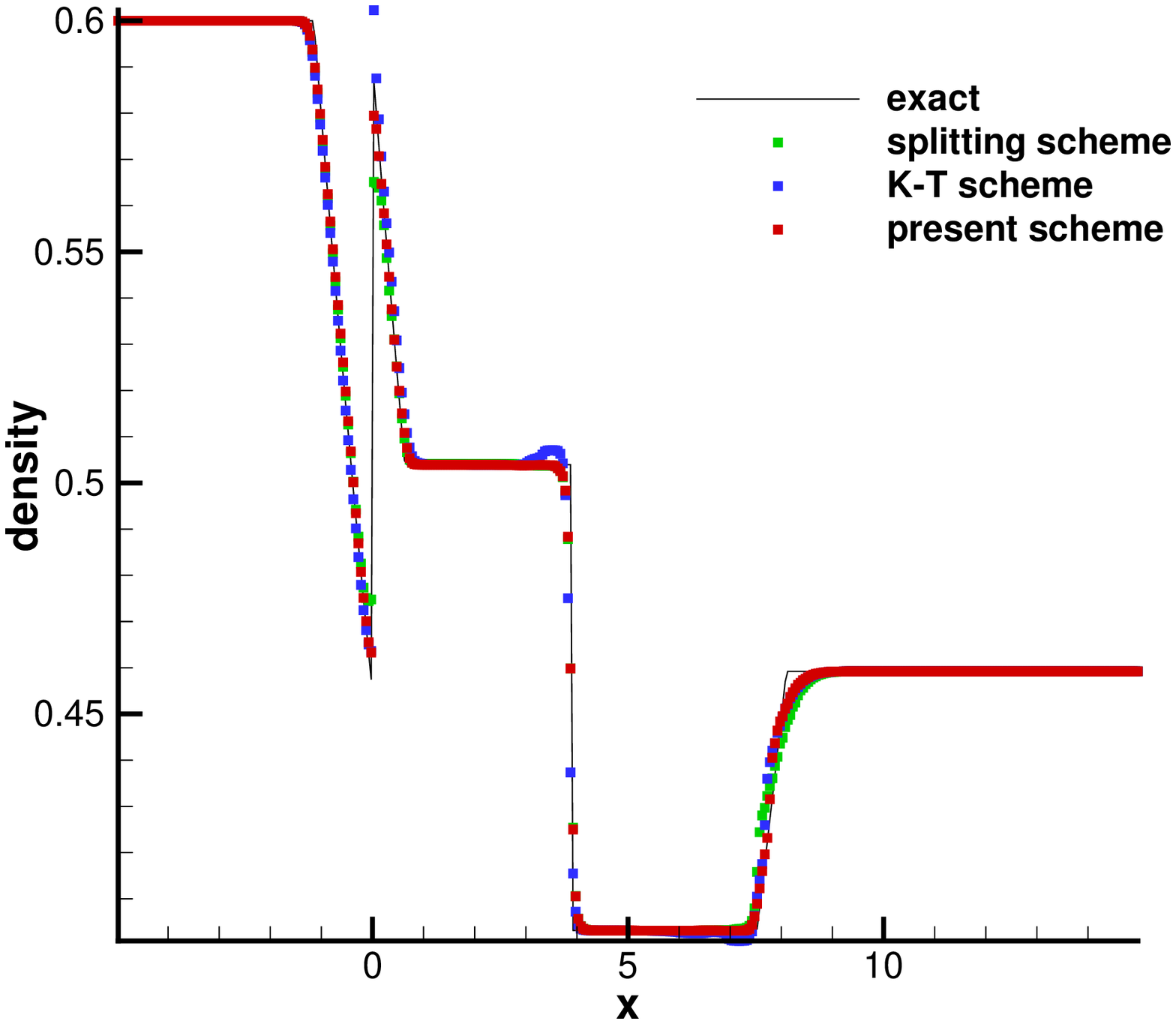}
	}
	\quad
	\subfigure[pressure]{
		\includegraphics[width=0.4\linewidth]{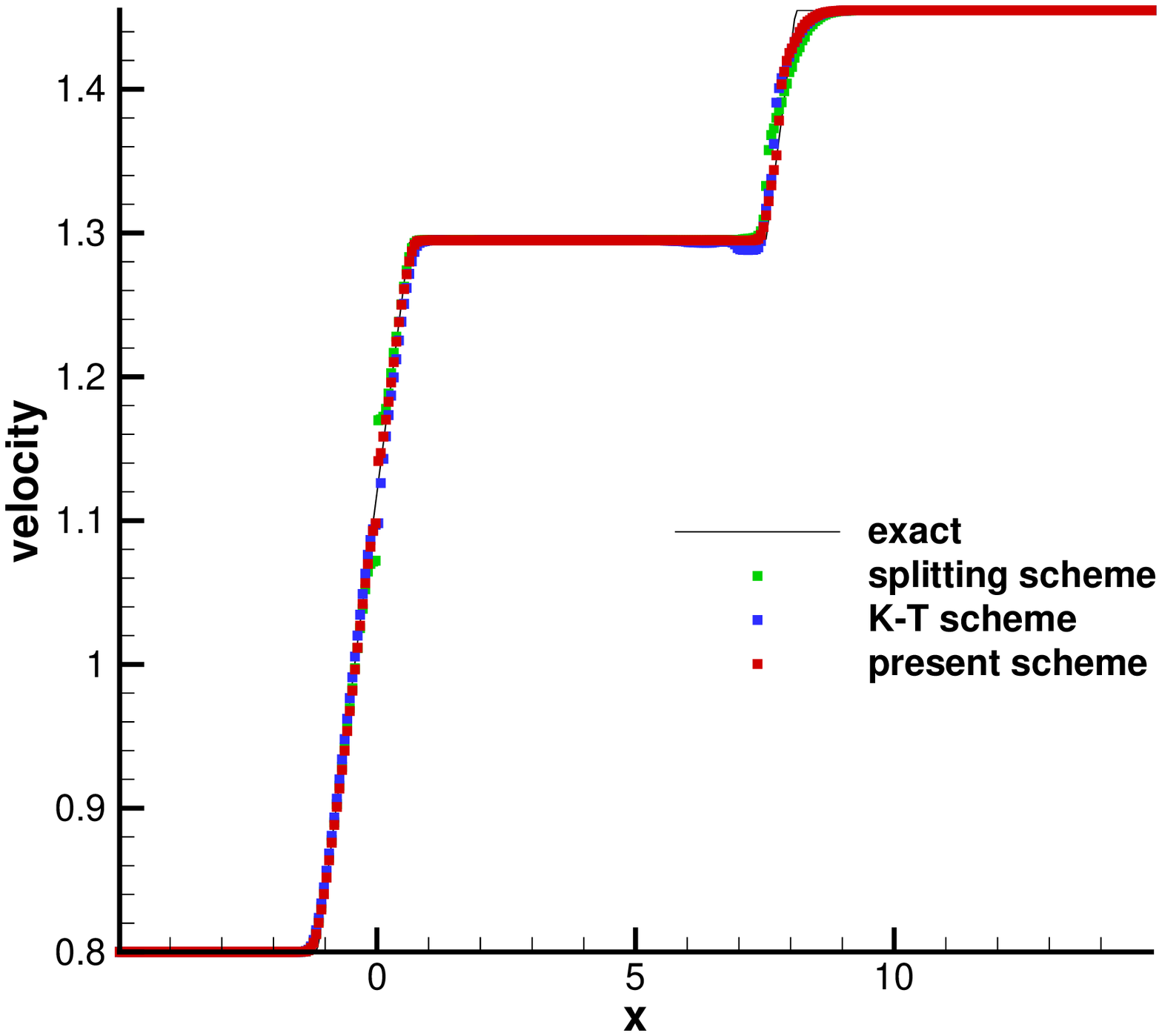}
	}
	\quad
	\subfigure[pressure]{
		\includegraphics[width=0.4\linewidth]{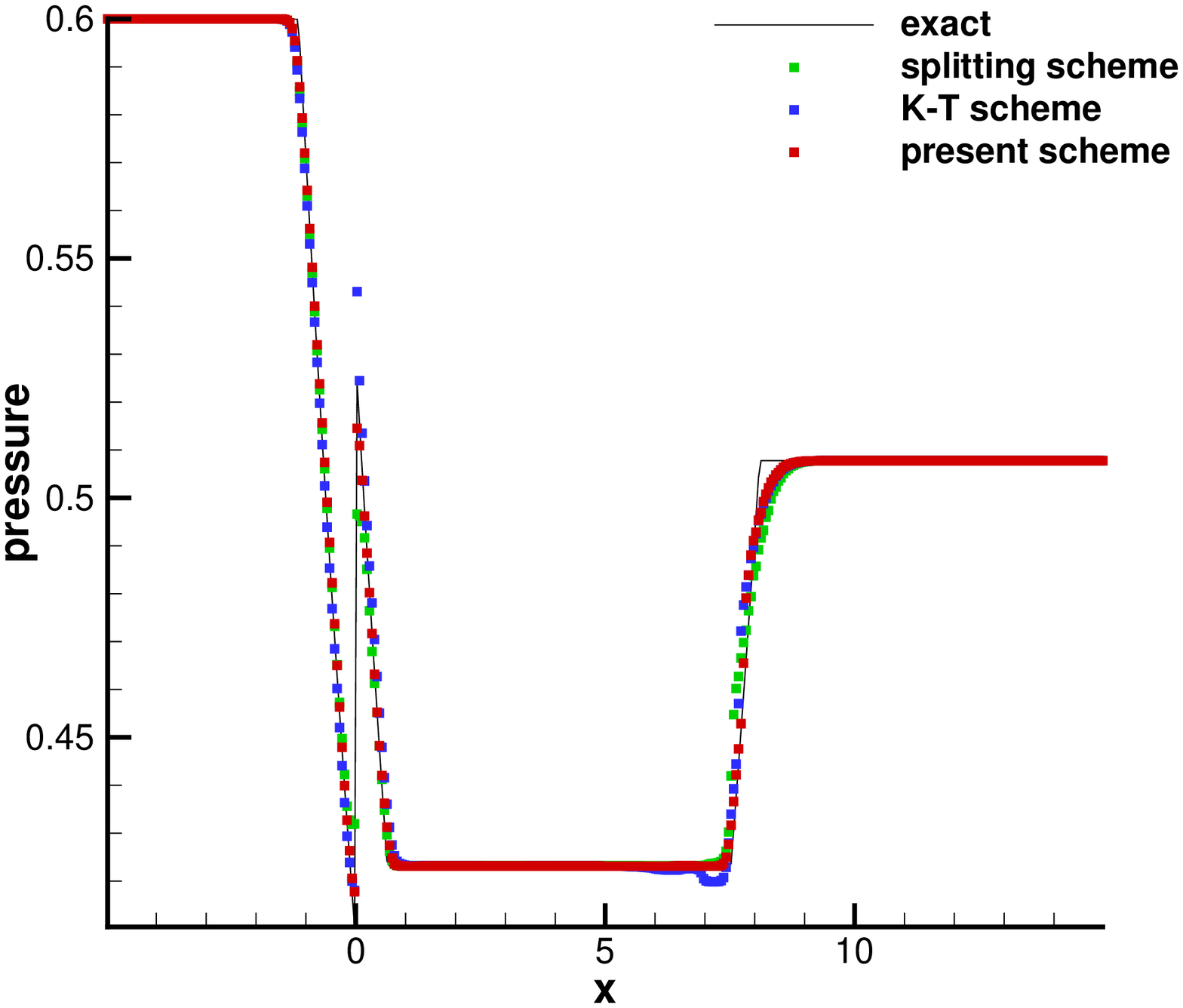}
	}
	\caption{Splitting scheme(green symbols), K-T scheme(blue symbols) and present scheme(red symbols) versus exact solutions(solid line) of Test8 with $h=0.05$ at time $t=3.0s$. Top left: density, Top right: velocity, and bottom: pressure.}
	\label{test8}
\end{figure}

The solutions of these tests are shown in Figure2, Figure\ref{test1-2}, Figure\ref{test1-3}, Figure\ref{test2-3}, Figure\ref{test4-5}, Figure\ref{test6-7}, Figure\ref{test8}. 

In Test1, the initial value is a single stationary wave. We use this test to check whether the numerical scheme is well-balanced. 
As shown in Figure\ref{test1-1}, the splitting scheme is not well-balanced. 
the incorrect decomposition of the convection equation for the initial discontinuity leads to the numerical solution containing several incorrect perturbations, which associates with three classical characteristic domains. 
In addition, the numerical solution has overshoots near the origin due to the inability to maintain steady solutions. 
Furthermore, neither of these errors is significantly reduced as the grid is refined. 
The K-T scheme and our sover-based scheme are both well-balanced (see Figure\ref{test1-2}, Figure\ref{test1-3}).

The exact solutions for both Test 2 and Test 3 are non-choked, with the former having a Type1 structure and the latter a Type2.
Figure\ref{test2-3} gives the results of the three numerical schemes for these two tests, all of which are in good agreement with the exact solution.
A disadvantage of the splitting scheme is that the numerical solution has overshoots in the vicinity of the stantionary wave.
The K-T scheme has a slight oscillation at the downstream contact discontinuity in the numerical solution for Test2. The results of the present solver-based numerical scheme are free of both errors.

The structure of the exact solution of Test4 is Type3, which is choked. 
The results of the three numerical schemes are shown in the left figure of Figure\ref{test4-5}.
The splitting scheme gives a relatively accurate approximation of the Riemann solution in each intermediate region, although its results have large numerical dissipation and overshoots near the stationary wave. 
For the K-T scheme, the numerical flux cannot satisfy the conditions in Theorem\ref{theorem: K-T scheme} in this test.
In order to make the K-T scheme work, we make some corrections, including the value of $I$ in (\ref{stationary wave curves1}) and (\ref{stationary wave curves2}) being set to zero if it is not solvable, and another branch of the stationary wave being adopted if the branch in (\ref{W-B: L-F flux}) is not available. 
These corrections are also used in a few of the later tests.
With these corrections, the K-T scheme works but there is still some non-physical oscillation.
The present solver-based scheme approximates each elementary wave better and has less numerical dissipation.

The exact solution for Test5 is a choked solution of Type4 structure and the right figure of Figure\ref{test4-5} shows the results for the three numerical schemes.
The solution of the splitting scheme converges to the exact solution, but with larger overshoots near the stationary wave.
The result of K-T scheme contains very large numerical oscillations, which we suspect, are due to the incorrect choice of branches of stationary wave curve.
The state of exact solution to the right of the origin is sonic, so the state of numerical solution to the right of the origin may be perturbed around the sonic state.
The perturbation of the numerical solution in the K-T scheme between the subsonic and supersonic states causes a change in the employed branch, which results in very large numerical oscillations.
The results of the present  solver-based scheme agrees well with the exact solution. Although our approximate Riemann solver ignores the Type4 structure, we achieve an approximation to the solution of this limit structure by the solution of the Type3 or Type2.

Test6 and Test7 are two choked solutions for $k<0$. 
The results of three schemes are shown in Figure\ref{test6-7}.
The splitting method contains greater numerical dissipation than the other two schemes, and the numerical solution in the constant region to the right of the origin in Test7 contains oscillations.
In both tests the K-T scheme gives incorrect results, and the results of our solver-based scheme agree well with the exact solution.

The last Test8 is an example of a choked solution for $k = 0$.
The exact solution has a rarefaction wave on each side of the origin in contact with it.
The numerical solutions for all three schemes are capable of approximating the exact solution, as shown in Figure\ref{test8}.
As can be seen from the density and pressure results, the present solver-based scheme works better than the other two schemes in the simulation of two rarefaction waves.

\section{Conclusions}\label{Conclusions}
For Euler equations with a singular source, a solution-structure based approximate Riemann solver was proposed in this paper.
The proposed solver is able to give the exact solution when the initial value is a single stationary wave and reasonable approximate solutions for different structures when the initial value is general. 
The numerical method based on our solver is not only well-balanced, but also can give numerical solutions in agreement with the exact solution in both non-extreme and extreme tests.

In one-dimensional problems we fixed the stationary singular source at a cell boundary. For one-dimensional equations with moving singular sources or for multidimensional equations, our solver can be directly generalised by means of reconstruction.

\section*{Acknowledgments}
This work was supported by the National Natural Science Foundation of China (Grant No.12101029) and Postdoctoral Science Foundation of China (Grant No.2020M680283).

\bibliographystyle{plain}
\bibliography{references}

\end{document}